\documentclass[a4paper,fleqn]{cas-sc}

\usepackage[numbers]{natbib}
\usepackage{amsmath}
\usepackage{kotex}
\usepackage{multirow}
\usepackage{graphicx}
\usepackage{amssymb}
\usepackage{relsize}
\usepackage{color}
\usepackage{amsfonts}
\usepackage[LGR,T1]{fontenc}
\usepackage{bigstrut}
\usepackage{array}
\usepackage{algorithm}
\usepackage{algpseudocode}
\usepackage{subcaption}
\usepackage{multicol}
\usepackage{mathtools}

\newcommand{\MR}[1]{\multirow{2}{*}{$#1$}}

\newcommand{\MCC}[1]{\multicolumn{3}{c}{#1}}
\newcommand{\MCCC}[1]{\multicolumn{4}{c}{#1}}
\newcommand{\mb}[1]{\mathbf{#1}}

\def\tsc#1{\csdef{#1}{\textsc{\lowercase{#1}}\xspace}}
\tsc{WGM}
\tsc{QE}
\tsc{EP}
\tsc{PMS}
\tsc{BEC}
\tsc{DE}

\begin{document}
\let\WriteBookmarks\relax
\def\floatpagepagefraction{1}
\def\textpagefraction{.001}
\shorttitle{Optimal hybrid parameter selection for stable sequential solution}
\shortauthors{Chang-uk Ahn et~al.}

\title [mode = title]{Optimal hybrid parameter selection for stable sequential solution of inverse heat conduction problem}

\author[1,2]{Chang-uk Ahn}
\address[1]{Department of Mechanical Engineering (Integrated Engineering), Kyung Hee University, 1732, Deogyeong-daero, Giheung-gu, Yongin-si, Gyeonggi-do 17104, Korea}
\ead{changuk.ahn@gmail.com}

\author[2]{Chanhun Park}
\ead{chpark@kimm.re.kr}

\author[2]{Dong Il Park}
\cormark[1]
\ead{parkstar@kimm.re.kr}
\address[2]{Department of Robotics $\And$ Mechatronics Research, Korea Institute of Machinery and Materials, Daejeon, 34103, Republic of Korea}

\author[1]{Jin-Gyun Kim}
\cormark[1]
\ead{jingyun.kim@khu.ac.kr}

\cortext[cor1]{Corresponding author}

\begin{abstract}
To deal with the ill-posed nature of the inverse heat conduction problem (IHCP), the regularization parameter $\alpha$ can be incorporated into a minimization problem, which is known as Tikhonov regularization method, popular technique to obtain stable sequential solutions.
Because $\alpha$ is a penalty term, its excessive use may cause large bias errors.
A ridge regression was developed as an estimator of the optimal $\alpha$ to minimize the magnitude of a gain coefficient matrix appropriately.
However, the sensitivity coefficient matrix included in the gain coefficient matrix depends on the time integrator; thus, certain parameters of the time integrators should be carefully considered with $\alpha$ to handle instability. 
Based on this motivation, we propose an effective iterative hybrid parameter selection algorithm to obtain stable inverse solutions.
We considered the Euler time integrator to solve IHCP using the finite element method.
We then considered $\beta$, a parameter to define Forward to Backward Euler time integrators, as a hybrid parameter with $\alpha$.
The error amplified by the inverse algorithm can be controlled by $\alpha$ first by assuming $\beta = 1$.
The total error is then classified into bias and variance errors.
The bias error can be computed using the maximum heat flux change, and the variance error can be calculated using the measurement noise error generated by prior information.
Therefore, $\alpha$ can initially be efficiently defined by the summation of the bias and variance errors computed in a time-independent manner.
Reducing the total error for better stability of the inverse solutions is also available by adjusting $\beta$, which is defined to minimize the magnitude of gain coefficient matrix when spectral radius of the amplification matrix is less than one.
Consequently, $\alpha$ could be updated with new $\beta$ in the iteration process.
The proposed efficient ridge estimator is essential to implement the iterative hybrid parameter selection algorithm in engineering practice.
The possibility and performance of the hybrid parameter selection algorithm were evaluated by well-constructed 1D and 2D numerical examples.
\end{abstract}

%


\maketitle

\section{Introduction} \label{Intro}

The inverse heat conduction problem (IHCP) is essential in most branches of science and engineering to estimate causal characteristics from measured temperatures.
Many analytical~\cite{Aviles1998, Monde2000, Cole2011} and numerical methods have been presented to solve the IHCP.
Numerical methods include three major discretization methods, namely, finite element method (FEM)~\cite{Tseng1995, Ling2003, Ling2005, Tseng1996, Krutz1978}, finite difference method (FDM)~\cite{Raynaud1986, Guo1991, Lin2004}, and boundary element method (BEM)~\cite{Kurpisz1992, Pasquetti1991, Singh2001}.
In addition, model reduction techniques~\cite{Glouannec2021, Lemonnier2021, NMPark2003, Tandy1986, Videcoq2001, Ahn2021} and parallel algorithm~\cite{Szenasi2018} have been developed in thermal sciences to handle large and complex discretized thermal systems.
In mathematics, the IHCP is called an ill-posed problem in the sense of Hadamard~\cite{Hadamard1902}. 
Inverse solutions are extremely sensitive to measurement errors, and thus small errors in the data are amplified dramatically~\cite{Tikhonov1997}.
To reduce the instability caused by the ill-posedness of inverse problems, various approaches have been developed for the stabilization of this discretized system, such as the sequential function specification method (SFSM)~\cite{Beck1985}, regularization method (RM)~\cite{Tikhonov1995, Tikhonov1997}, iterative regularization~\cite{Alifanov1994, Alifanov1989, Rabczuk2017, Reinhart1996}, combined SFSM--RM~\cite{Beck1986}, Levenberg-Marquardt method~\cite{Ozisik2000}, repulsive particle swarm optimization (RPSO) method~\cite{KHLee2008, KHLee2019}, and stochastic method~\cite{Narayanan2004}.

The inverse problems are classified into whole time-domain and sequential time-domain approaches~\cite{Ozisik2000, Woodbury1998, Woodbury2002}.
The whole-time domain method is relatively more stable than the sequential time-domain method because it uses the sum of the data in which errors are smoothed.
Because the whole time-domain method is based on data averaging approaches, it is effective for the parameter estimation and also provides good stability.
On the other hand, the sequential time-domain method is an effective way to find time-dependent unknown functions, such as temperature and heat flux.
In the sequential-time domain method, the unknown time-dependent response can be represented as a polynomial in time with a finite number of certain parameters.
Therefore, it is useful for estimating the time-varying temperature and heat flux that are difficult to measure.
However, it allows for an exact match between the estimated and the measured temperature at every time step, causing severe oscillation in inverse solutions.
The oscillation behavior becomes more unstable when high-frequency error data with a small time interval are used for good resolution.

Tikhonov regularization (TR) method~\cite{Tikhonov1995, Tikhonov1997}, one of the most popular methods, is simple and effective in obtaining stable sequential solutions.
The iterative regularization (IR) method proposed by Alifanov~\cite{Alifanov1994, Alifanov1989} is also a powerful approach.
Still, TR method has a lower computational cost than IR method if only the regularization parameter can be obtained in advance; thus, it could be utilized in real-time measurement or control environments.
TR method involves adding penalty terms to a minimization problem to reduce excessive sensitivity to measurement noise, resulting in smooth and stable solutions.
An issue exists that requires proper selection of the optimal regularization parameter $\alpha$.
This is because $\alpha$ balances the trade-off between filtering out sufficient variance error and not losing significant amount of information in the solution, which means an increase in the bias error.
To date, general techniques are available for selecting the optimal $\alpha$, such as the discrepancy principle~\cite{Groetsch1983, Engl1987, Bauer2011}, generalized cross-validation (GCV)~\cite{Golub1979, Trujillo1989, Lukas2006}, and L-curve rule~\cite{Hansen1993, Trujillo1994}.
The optimal $\alpha$ estimators are also called the ridge estimators.

Optimal $\alpha$ is the value obtained when the magnitude of a gain coefficient matrix is appropriately minimized.
The gain coefficient matrix determines the stability and accuracy in inverse problems~\cite{Beck1985}.
However, a sensitivity coefficient matrix included in the gain coefficient matrix depends on the time integrator: thus, certain parameters of the time integrators should be carefully considered with $\alpha$ to handle instability. 
In this work, the Euler time integrator is used, and thus its parameter $\beta$ should be considered with $\alpha$.
The parameter $\beta$ allows selection from the Forward Euler ($\beta=0$) to the Backward Euler ($\beta=1$) method. 
Further, $\beta$ can be selected from the range of 0 to 1 so that the magnitude of the gain coefficient matrix is minimal; however, the Backward Euler method with $\beta = 1$ has been used in general to avoid divergence of the solution, although accuracy is lost.
Only a few studies have presented the general stability conditions~\cite{Liu1996, Reinhardt1994, Ling2006}.

Motivated by these factors, an iterative hybrid parameter ($\alpha$ and $\beta$) selection algorithm is presented in this study to efficiently obtain stable solutions.
A fast ridge estimator to effectively select the optimal $\alpha$ is a prerequisite to accomplish this objective.
Additionally, a new stability condition of the inverse problem with the Tikhonov regularization is required to select $\beta$. 
The discrepancy principle is one of the most popular ridge estimators because it is the simplest method that also yields good accuracy. 
However, the discrepancy principle needs to solve a direct problem at each $\alpha$ discretized over a specific range; thus, it is not effective for use with the proposed iterative hybrid parameter selection algorithm.
Further, the discrepancy principle is typically based on the measured temperature data, which may not guarantee the result of heat flux estimation.

To overcome this limitation, we propose a new ridge estimator for heat flux---it starts with decomposing the total error in the heat flux into bias and variance errors.
Employing standard statistical assumptions indicates that the bias and variance errors respectively depend on heat flux change and measurement noise, and the relationship between these two errors is linearly independent.
Beck suggested this idea to find the optimal $\alpha$ effectively in a numerical setting~\cite{Woodbury2013}.
The ridge estimator proposed by Beck is based on the residual sum of squares in the errors in the heat flux, which is the best measure to recover the heat flux history as accurately as possible.
However, the level of error in the heat flux is rarely known, and the exact heat flux history is generally unknown.
Accordingly, using Beck's method is perhaps difficult to evaluate the error measure of the heat flux~\cite{Woodbury2013}.
The proposed estimator is also based on the expected error of the heat flux, but the maximum heat flux change and the temperature deviation are used unlike the Beck's method.
The heat flux deviation results from the amplification of the measurement noise of the temperature in the inverse problem and, thus, can be calculated by the temperature deviation, which is generally known.
Therefore, the measure of the proposed estimator is the sum of the bias and variance errors, which are calculated by the maximum heat flux change and the temperature deviation, respectively.
In detail, only maximum bias error is computed with the maximum heat flux change, instead of computing the bias error propagation over time.
This is because the bias error depends significantly on heat flux changes in a well-established experimental environment, with less influence by ill-conditioning of the inverse solver.
The variance error propagation can be approximately represented as a form of power series in the measurement noise of the temperature.
This allows us to calculate the propagated variance error only by matrix multiplication without time integration.
Therefore, the proposed ridge estimator could be used in practice in a time-independent manner.

To select $\beta$, we introduce the new stability condition for $\beta$ based on spectral norm analysis.
When $\beta < 0.5$, it is not unconditionally stable.
Therefore, the time step size, mesh size, sensor numbers, and sensor locations could be adjusted for stability.
These conditions are associated with hardware equipments and numerical modeling, so it is not easy to control in inverse problems.
For these reasons, the stability condition adjusted by $\beta$ is much more efficient in practice.
To do this, in inverse problems, input errors (i.e., noise data measured sequentially at specific time intervals), as well as initial error propagation, should be considered for stability.
Because it is challenging to construct an amplification matrix that simultaneously represents the propagation of these two errors, we control the stability of the input and initial errors separately.
Therefore, the amplification matrix for the initial error is newly derived here, and then $\beta$ is selected within the range satisfying the condition that the spectral radius of the amplification matrix is less than 1~\cite{Hirsch1988}. 
The instability for input errors is handled by the proposed ridge estimator.
Finally, in the proposed iterative algorithm, the optimal $\alpha$ is defined first by the efficient ridge estimator by assuming $\beta = 1$.
Then $\beta$ is determined by considering the magnitude of the gain coefficient matrix and the proposed stability condition, 
and then insufficient stability is regulated by $\alpha$ with the new $\beta$.
The algorithm is repeated if necessary.
The results demonstrated that the proposed ridge estimator is more effective than conventional ridge estimators; therefore, the iterative hybrid parameter selection algorithm is useful for estimating accurate and stable sequential solutions.

The remainder of this paper is organized as follows. 
Section ~\ref{DirectInverse} presents the direct and inverse problem based on FEM.
Section ~\ref{hybrid} introduces the algorithm of the hybrid parameter selection.
In detail, the efficient ridge estimator is proposed to estimate the optimal $\alpha$ in Section~\ref{RegularAlpha}.
The new amplification matrix of the inverse problem is derived for the stability condition to select the parameter $\beta$ in Section~\ref{Beta}.
The iterative hybrid parameter selection algorithm is organized in Section~\ref{Algorithm0}.
Section~\ref{NumericalExamples} evaluates the feasibility and performance of the proposed algorithm by well-constructed 1D and 2D numerical examples.
Finally, Section~\ref{Conclusion} presents the conclusions.

\section{Direct and Inverse Poblems}  \label{DirectInverse}
\subsection{Direct Poblem}  \label{DirectProblem}
\begin{figure}
	\centering
		\includegraphics[scale=.12]{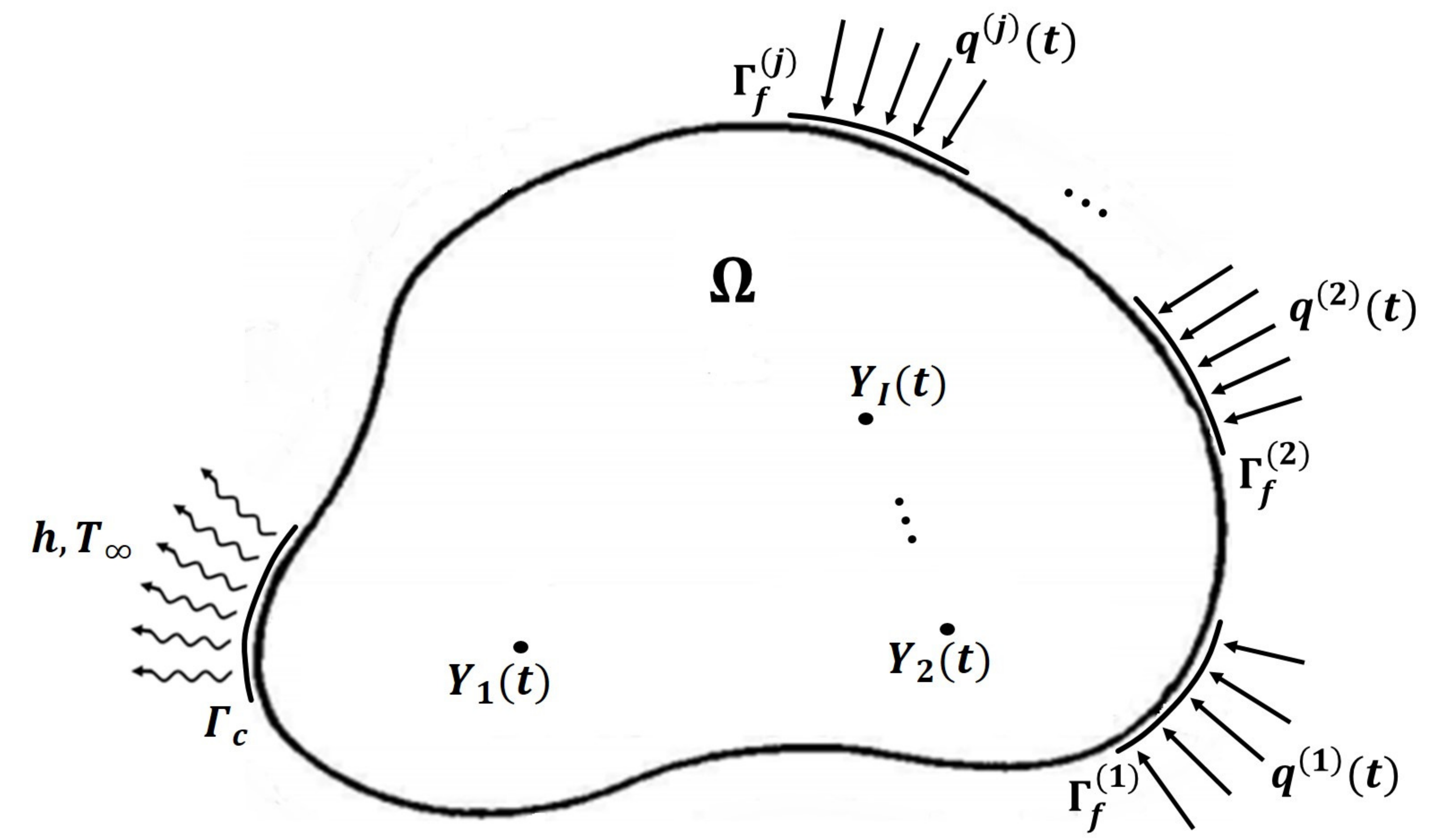}
	\caption{Schematic of the heat conduction domain with the measured temperature $Y_i(t)$, ($i=1,\cdots,I$).}
	\label{Domain}
\end{figure}

As a preliminary, we explain the FEM for the direct heat conduction problem.
The governing equations for the temperature fields, as given by the law of conservation of energy, consist of the following partial differential equations~\cite{Bathe1995}:
     \begin{equation} \label{Govern}
     \begin{cases} \rho c_p \frac{\partial T(x,y,z,t)}{\partial t} = \frac{\partial}{\partial x}\left[k_x\frac{\partial T(x,y,z,t)}{\partial x}\right] 
     + \frac{\partial}{\partial y}\left[k_y\frac{\partial T(x,y,z,t)}{\partial y}\right] + \frac{\partial}{\partial z}\left[k_z\frac{\partial T(x,y,z,t)}{\partial z}\right], \qquad & \mbox{in}~\Omega, \\[1.5ex]
     k \frac{\partial T(x,y,z,t)}{\partial \mb{n}} = -h\left[T(x,y,z,t)-T_{\infty}\right], & \mbox{on}~\Gamma_c, \\[1.5ex]
     k \frac{\partial T(x,y,z,t)}{\partial \mb{n}} = q^{(j)}(t), \qquad (j=1,2,\cdots,N), & \mbox{on}~\Gamma_{f}^{(j)}, \\[1.5ex]
     T(x,y,z,0) = T_0, & \mbox{in}~\Omega,
     \end{cases}
     \end{equation} \\
where $\rho$, $c_p$, $k_x$, $h$, $T_0$, and $T_\infty$ denote the mass density, specific heat capacity, thermal conductivity in the $x$ direction, convective heat transfer coefficient, initial temperature, and ambient temperature, respectively.
We assume homogeneous material properties, i.e., $k=k_x=k_y=k_z$.
The spatial domain is denoted by $\Omega$.
A heat flux profile $q^{(j)}(t)$ is specified on the $j$-th boundary, denoted by $\Gamma_{f}^{(j)}$, with an outward unit normal to $\Omega$ denoted by $\mb{n}$. 
A schematic diagram of the heat conduction problem is illustrated in Figure~\ref{Domain}.
Convection is specified on a boundary, denoted by $\Gamma_c$.

The partial differential equation is solved using the Galerkin finite element method.
The standard Galerkin procedure transforms as shown in the following algebraic differential equations:
     \begin{equation} \label{FEM}
     \mb{C}\mb{\dot{T}}(t)+\mb{K}\mb{T}(t)=\mb{q}_g(t),
     \end{equation}
with the components of the matrices and vector, which are expressed as
     \begin{subequations} \label{FEMcomp}
     \begin{gather} 
     \label{FEMcompC}
     \mb{C} = \sum_{e=1}^{N_e} \left[\int_{\Omega} \rho c_p \mb{L}^{(e)^T} \mb{H}^T \mb{H} \mb{L}^{(e)} dV\right],\\
     \label{FEMcompK}
     \mb{K} = \sum_{e=1}^{N_e} \left[\int_{\Omega} \mb{L}^{(e)^T} \mb{B}^T \mb{D B} \mb{L}^{(e)} dV \right] + \sum_{e=1}^{N_c} \left[\int_{\Gamma_c} h \mb{L}_c^{(e)^T} \mb{H}^T \mb{H} \mb{L}_c^{(e)} dS \right], \\
     \label{FEMcompq}
     \mb{q}_g(t) = \sum_{j=1}^{N} \sum_{e=1}^{N_{f}^{(j)}} \left[\int_{\Gamma_{f}^{(j)}} q^{(j)}(t) \mb{L}_{f}^{(e,j)^T} \mb{H}^T dS\right] + \sum_{e=1}^{N_c} \left[\int_{\Gamma_c} h \mb{L}_c^{(e)^T} \mb{H}^T \mb{T}_{\infty} dS \right],
     \end{gather}
     \end{subequations}
where $\mb{C}$, $\mb{K}$, $\mb{H}$, $\mb{B}$, and $\mb{D}$ denote the global heat capacity, global heat transfer, 
shape function, spatial derivative of the shape function, and material property matrices, respectively. 
$\mb{T}(t)$ and $\mb{q}_g(t)$ denote the global nodal temperature and heat flux vector, respectively.
In the above formulation, the global heat transfer matrix can be found as the sum of the global conductivity and convection matrices.
$\mb{L}^{(e)}$ is a Boolean matrix establishing the relation between the element $(e)$ and the assembled (global) degrees of freedom.
$\mb{L}_{c}^{(e)}$ and $\mb{L}_{f}^{(e,j)}$ are built for mapping to the boundary of the convection and heat flux, respectively.
$N_e$ and $N_g$ are the total numbers of elements and degrees of freedoms (DOFs).
The numbers of element specified on the boundary of the convection and heat flux are respectively defined as $N_c$ and $N_f^{(j)}$.
The local index $(j)$ spans from 1 to $N$. The local index $(e)$ spans from 1 to $N_e$, $N_c$, or $N_f^{(j)}$.

To solve Eq.~\ref{FEM} Euler method can be employed, which is given by
     \begin{equation} \label{TimeIntegrator}
     \mb{T}^{m+1} = \mb{T}^m + \Delta t \left[ \left( 1 - \beta \right) \mb{\dot{T}}^m + \beta \mb{\dot{T}}^{m+1}\right],  \qquad 0\le\beta\le1.
     \end{equation} 

\begin{figure}
	\centering
		\includegraphics[scale=.3]{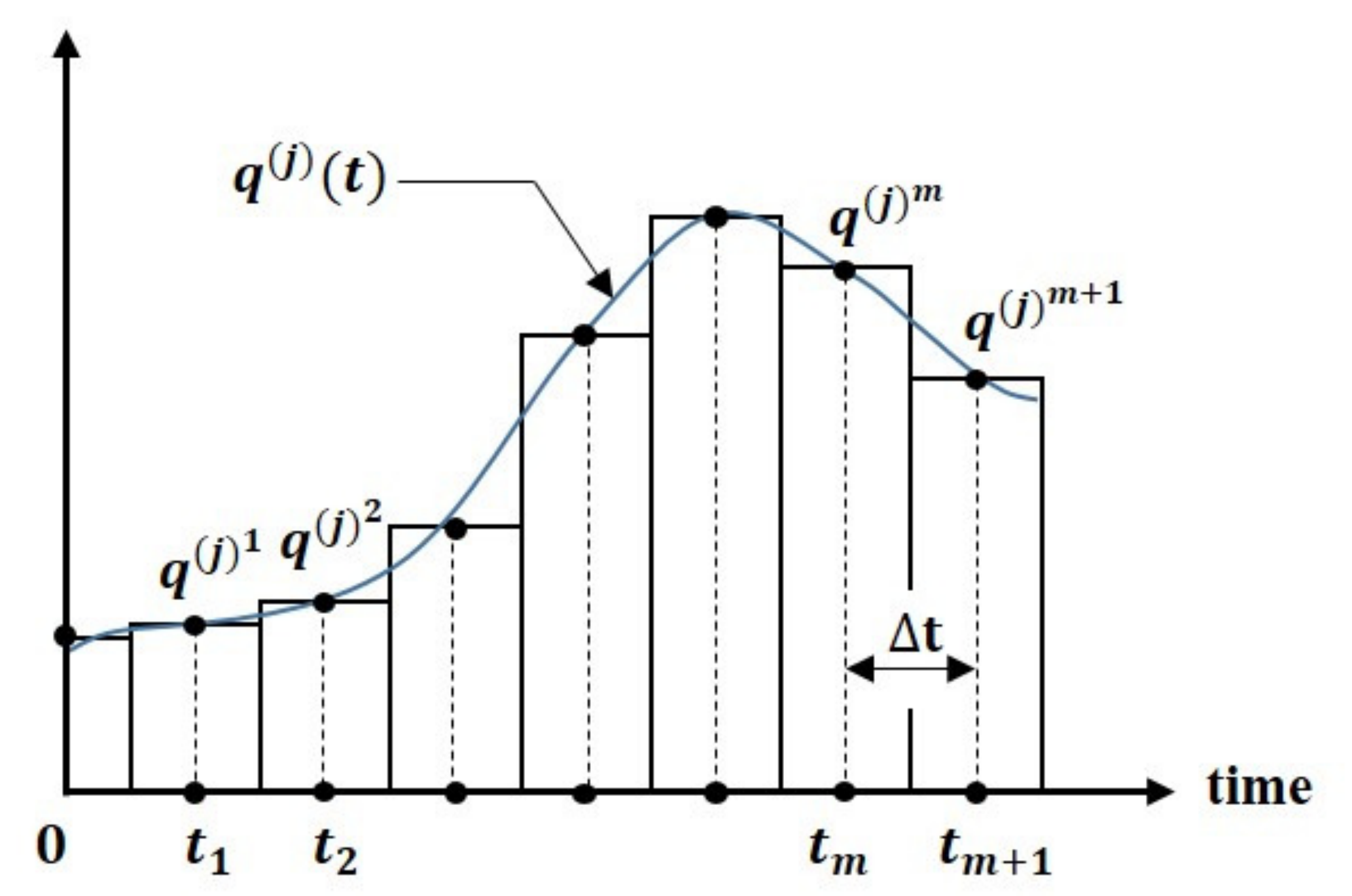}
	\caption{Constant piecewise functional form of the heat flux, which is specified on the $j$-th boundary and defined by the trapezoidal rule ($\beta=0.5$).}
	\label{TimeStep}
\end{figure}

Then, we obtain
     \begin{equation} \label{Euler}
     \left(\frac{1}{\Delta t}\mb{C} + \beta\mb{K} \right)\mb{T}^{m+1} = \left[\frac{1}{\Delta t}\mb{C} - \left(1-\beta\right)\mb{K}\right]\mb{T}^{m} + \left[ \left(1 - \beta \right)\mb{q}_g^m + \beta\mb{q}_g^{m+1}\right],
     \end{equation} 
where the time $t$ is represented in a discretized form $t_m$, and $\Delta t$ denotes the time--step size, which is defined as $(\Delta t = t_{m+1}-t_m)$.
In general, the time--step size is set equivalent to a measurement time interval (or sampling time interval).
The parameter $\beta$ determines the performance of the time integrator, which is the Backward Euler method when $\beta=1$ and Forward Euler method when $\beta=0$.
When using a lumped capacitance model, caution must be exercised if $\beta=0$.
An unknown heat flux history can be estimated if the surface on which the heat fluxes are specified contains the subsurface where the temperatures are measured.
Otherwise, the estimated heat flux would be zero.
Equation~\ref{Euler} can be rewritten as a standard temperature equation in $\mb{T}^{m+1}$:
     \begin{subequations} \label{Tg}
     \begin{gather} 
     \label{Tgm}
     \mb{T}^{m+1}  = \mb{A}\mb{T}^{m} + \mb{U}\left[ \left(1 - \beta \right)\mb{q}_g^m + \beta\mb{q}_g^{m+1}\right],\\
     \label{TgA}
     \mb{A} = \mb{U}\left[\frac{1}{\Delta t}\mb{C} - \left(1-\beta\right)\mb{K}\right],\\
     \label{TgU}
     \mb{U}  = \left(\frac{1}{\Delta t}\mb{C} + \beta\mb{K} \right)^{-1},
     \end{gather}
     \end{subequations}
where the heat flux profile, $q^{(j)}(t)$ in Eq.~\ref{FEMcompq}, is approximated to a constant piecewise functional form $q^{(j)^m}$ as shown in Figure~\ref{TimeStep}.
In other words, the global nodal heat flux vector $\mb{q}_g^{m}$ is constructed by $\mb{q}^m$, which is represented as follows:
     \begin{equation} \label{HeatFlux}
     \mb{q}^m = 
     \begin{bmatrix}
     {q}^{(1)^m}    &  {q}^{(2)^m} &  \cdots   &  {q}^{(N)^m}
     \end{bmatrix}^T.
     \end{equation}

The matrix normal equation and the unknown heat flux component $\mb{q}^{m+1}$ could be obtained based on the standard temperature equation in Eq.~\ref{Tg}.

\subsection{Inverse Problem with Regularization Method}  \label{ClassicMethod}
In this section, we review the sequential method for the IHCP, which estimates an unknown heat flux components at every time interval.
The inverse problems are solved by minimizing the sum of square error, which is given by
     \begin{equation} \label{least}
     S^{m+1}=\left[\mb{Y}^{m+1}-\mb{T}^{m+1}\right]^T\left[\mb{Y}^{m+1}-\mb{T}^{m+1}\right], \quad \mb{Y}^{m+1} = [Y_1^{m+1}, \cdots, Y_I^{m+1}]^T, \quad \mb{T}^{m+1} = [T_1^{m+1}, \cdots, T_I^{m+1}]^T,
     \end{equation}
where $S^{m+1}$ is an objective function. $\mb{Y}^{m+1}$ and $\mb{T}^{m+1}$ are the vector of measured temperature and the corresponding set of the vector of calculated temperature, respectively.
$I$ is the total numbers of sensors.

The minimization of $S^{m+1}$ can be represented by setting the gradient of $S^{m+1}$ to zero:
     \begin{equation}  \label{GradS}
     \nabla \mb{S}^{m+1} = -2\mb{X}^{m+1^T}\left[\mb{Y}^{m+1}-\mb{T}^{m+1}\right] = \mb{0},
     \qquad
     \mb{X}^{m+1} = \frac{\partial \mb{T}^{m+1}}{\partial \mb{q}^{m+1}} = 
     \begin{bmatrix}
     \frac{\partial \mb{T}^{m+1}}{\partial q^{(1)^{m+1}}}, & \cdots, & \frac{\partial \mb{T}^{m+1}}{\partial q^{(N)^{m+1}}}
     \end{bmatrix},
     \end{equation}
where $\mb{X}^{m+1} $ is a sensitivity (or Jacobian) coefficient matrix and is defined as the first derivative of $\mb{T}^{m+1} $ with respect to $\mb{q}^{m+1}$.
The matrix size of $\mb{X}^{m+1}$ is $I\times N_g$. The vector size of $\mb{Y}^{m+1}$ and $\mb{T}^{m+1}$ are $I\times 1$ and $N_g\times 1$, respectively.
This indicates $I = N_g$, which means the number of sensors is the same as the number of total DOFs.  
However, the number of sensors is much smaller than the number of total DOFs because  the sensors are attached to a few locations in practice.
Therefore, Eq.~\ref{GradS} should be redefined as
     \begin{equation} \label{GradSsDetail} 
     \nabla \mb{S}^{m+1} = -2\mb{X}_s^{m+1^{T}}\left[\mb{Y}^{m+1} - \mb{L}_s\mb{T}^{m+1}\right] = \mb{0}, \qquad \mb{X}_s^{m+1} = \mb{L}_s\mb{X}^{m+1},
     \end{equation}
where $\mb{L}_s$ is the $I\times N_g$ Boolean matrix, which leaves only the elements of $\mb{T}^{m+1}$ and $\mb{X}^{m+1}$ corresponding to the sensor positions.

To determine the unknown heat flux component, $\mb{T}^{m+1}$ is expressed by the Taylor expansion, as follows.
     \begin{equation} \label{Taylor}
     \mb{T}^{m+1} \approx \bar{\mb{T}}^{m+1} + \mb{X}^{m+1} \Delta \mb{q}^{m+1},   \qquad \Delta \mb{q}^{m+1} = \mb{q}^{m+1} - \mb{q}^{m}. 
     \end{equation}

Considering Eq.~\ref{Tg}, the components are given as
     \begin{subequations} \label{TaylorSet}
     \begin{gather} 
     \label{TaylorSetT}
     \bar{\mb{T}}^{m+1} =\mb{T}^{m+1}\mid_{\mb{q}^{m+1}=\mb{q}^m} = \mb{A}\mb{T}^{m} + \mb{U}\mb{q}_g^{m}, \\
     \label{TayloSetX}
     \mb{X}^{m+1} = \frac{\partial \mb{T}^{m+1}}{\partial \mb{q}^{m+1}}
                      = \mb{U} \frac{\partial \mb{q}_g^{m+1}}{\partial \mb{q}^{m+1}},\\
     \label{UnitQ}
     \frac{\partial \mb{q}_g^{m+1}}{\partial \mb{q}^{m+1}} = 
     \left[ \sum_{e=1}^{N_{f}^{(1)}} \left[\int_{\Gamma_f^{(j)}} \mb{L}_{f}^{(e,1)^T} \mb{H}^T dS \right],  \cdots, \sum_{e=1}^{N_{f}^{(N)}} \left[\int_{\Gamma_f^{(j)}} \mb{L}_{f}^{(e,N)^T} \mb{H}^T dS \right] \right],
     \end{gather}
     \end{subequations}
where $\bar{\mb{T}}^{m+1}$ is the virtual temperature at a uniform load condition ($\mb{q}^{m+1}=\mb{q}^m$), 
and $\Delta \mb{q}^{m+1}$ is the amount of change in the heat flux components from $\mb{q}^m$ to $\mb{q}^{m+1}$.
Here, $\left(\partial \mb{q}_g^{m+1}/\partial \mb{q}^{m+1}\right)$ is a unit heat flux matrix.
In this study, we focus on the linear problem, i.e., the thermal material properties are constants.
Therefore, $\mb{A}$ and $\mb{U}$ in Eq.~\ref{TaylorSet} do not change with time; thus, $\mb{X}^{m+1}$ and $\left(\partial \mb{q}_g^{m+1}/\partial \mb{q}^{m+1}\right)$ are denoted as $\mb{X}$ and $\left(\partial \mb{q}_g/\partial \mb{q}\right)$, respectively.

Substituting Eq.~\ref{Taylor} into Eq.~\ref{GradSsDetail} gives
     \begin{equation} \label{MNE}
     \Delta \mb{q}^{m+1} = \mb{G} \left[\mb{Y}^{m+1} - \mb{L}_s\bar{\mb{T}}^{m+1}\right], \qquad \mb{G} = \left(\mb{X}_s^T\mb{X}_s\right)^{-1}\mb{X}_s^T. 
     \end{equation}
This equation is known as a matrix normal equation, and $\mb{G}$ is called a gain coefficient matrix~\cite{Beck1985}. 
Finally, the unknown vectors $\mb{q}^{m+1}$ and $\mb{T}^{m+1}$ can be estimated as follows.
     \begin{subequations} \label{EstTq}
     \begin{gather} 
     \label{Estq}
     \mb{q}^{m+1} = \mb{q}^{m} + \Delta \mb{q}^{m+1},\\
     \label{EstT}
     \mb{T}^{m+1} = \bar{\mb{T}}^{m+1} + \mb{X} \Delta \mb{q}^{m+1}.
     \end{gather}
     \end{subequations}

To obtain $\mb{q}^{m+1}$ and $\mb{T}^{m+1}$ in Eq.~\ref{EstTq}, $\mb{X}_s^T\mb{X}_s$ in Eq.~\ref{MNE} must be non-singular.
However, if $\mb{X}_s^T\mb{X}_s$ is close to the singular matrix, the solution $\mb{q}^{m+1}$ becomes unstable.
This problem can be alleviated by adding regularization terms to the objective function in Eq.~\ref{least}, as shown in the following representation~\cite{Tikhonov1995, Tikhonov1997}:
     \begin{equation} \label{leastRegular}
     S_{\alpha}^{m+1}=\left(\mb{Y}^{m+1}-\mb{L}_s\mb{T}_{\alpha}^{m+1}\right)^T\left(\mb{Y}^{m+1}-\mb{L}_s\mb{T}_{\alpha}^{m+1}\right) 
             + \alpha \sum_{r=0}^{2}\left[W_r\left(\mb{H}_r\Delta\mb{q}_{\alpha}^{m+1}\right)^T\left(\mb{H}_r\Delta\mb{q}_{\alpha}^{m+1}\right)\right], 
     \end{equation}
where $\alpha$ is a regularization parameter. $S_{\alpha}^{m+1}$, $\mb{T}_{\alpha}^{m+1}$, and $\Delta\mb{q}_{\alpha}^{m+1}$ are functions of the regularization parameter.
$W_r$ and $\mb{H}_r$ are the weight factors and square matrices for the $r$-th order regularization method, respectively.
This study only introduces the zeroth-order ($r=0$) regularization method, where $\mb{H}_0$ is a unit matrix.
In a similar manner as before, the matrix normal equation and the estimated solutions are derived as
     \begin{subequations} \label{MNE_reg}
     \begin{gather}
     \label{MNEq_r}
     \Delta \mb{q}_{\alpha}^{m+1} = \mb{G}_{\alpha} \left[\mb{Y}^{m+1} - \mb{L}_s\bar{\mb{T}}_{\alpha}^{m+1}\right],  
     \qquad \mb{G}_{\alpha} = \left(\mb{X}_s^T\mb{X}_s + \alpha \mb{I}\right)^{-1}\mb{X}_s^T,\\
     \mb{q}_{\alpha}^{m+1} = \mb{q}_{\alpha}^{m} + \Delta \mb{q}_{\alpha}^{m+1},\\
     \label{Test}
     \mb{T}_{\alpha}^{m+1} = \bar{\mb{T}}_{\alpha}^{m+1} + \mb{X} \Delta \mb{q}_{\alpha}^{m+1},
     \end{gather}
     \end{subequations}
where the condition number of $\mb{X}_s^T\mb{X}_s$ is reduced by adding positive elements to its diagonals, following which the oscillatory behavior of the solution is damped.
For a large value of $\alpha$, the second term in Eq.~\ref{leastRegular} is dominant.
As a result, for $r=0$, $\Delta\mb{q}_{\alpha}^{m+1}$ tends to become constant, and the solution $\mb{q}_{\alpha}^{m+1}$ deviates from the exact solution \cite{Beck1985}.
Therefore, it is critical to develop the ridge estimator \cite{Golub1979}.
The Morozov discrepancy principle is a conventional approach to estimate the optimal $\alpha$, and is given by~\cite{Groetsch1983},
     \begin{equation} \label{Morozov}
     \frac{1}{M}\sum_{m=1}^{M}\left(Y^{m}-\mb{L}_{s}\mb{T}_{\alpha}^{m}\right)^T\left(Y^{m}-\mb{L}_{s}\mb{T}_{\alpha}^{m}\right) = \sigma^2,
     \end{equation}
where $\sigma$ is the priori known noise level in the measurement data, and
$M$ is an end--time step, i.e., $t_M$ is the total measurement time.
In this principle, the optimal $\alpha$ can be selected using the equality of Eq.~\ref{Morozov}.
This estimator is one of the simplest and most widely used strategies.
However, this method involves high computational cost because it requires full-time simulations at each $\alpha$ discretized over a specific range.
Section~\ref{RegularAlpha} proposes a method to alleviate this problem.

\section{Hybrid Parameter for Inverse Problem} \label{hybrid}
As discussed in Section~\ref{DirectInverse}, solving the inverse problem requires the time integrator and the objective function.
In this study, the Euler method and ordinary least square error with zeroth-order regularization are used:
     \begin{subequations} \label{IHCPfinal}
     \begin{gather} 
     \label{Solver}
     \textbf{Objective function}: S_{\alpha}^{m+1}=\left(\mb{Y}^{m+1}-\mb{L}_s\mb{T}_{\alpha}^{m+1}\right)^T\left(\mb{Y}^{m+1}-\mb{L}_s\mb{T}_{\alpha}^{m+1}\right) 
             + \alpha \left(\Delta\mb{q}_{\alpha}^{m+1}\right)^T\left(\Delta\mb{q}_{\alpha}^{m+1}\right),\\
     \label{TimeInt}
     \textbf{Time integrator}: \mb{T}^{m+1} = \mb{T}^m + \Delta t \left[ \left( 1 - \beta \right) \mb{\dot{T}}^m + \beta \mb{\dot{T}}^{m+1}\right].
     \end{gather}
     \end{subequations} 

The IHCP society has focused on finding the optimal $\alpha$ to reduce the variance error, and $\beta$ has been selected independently in the Backward Euler method ($\beta=1$) in general.
However, the stability and accuracy of the inverse solution depend highly on the combination of selecting $\alpha$ and $\beta$
because the gain coefficient matrix is a function of $\alpha$ and $\beta$.
Therefore, the relationship between $\alpha$ and $\beta$ should be considered in the inverse process.
The proposed iterative selection algorithm is based on this motivation. 
Section~\ref{RegularAlpha} proposes an efficient method to estimate the optimal $\alpha$.
Section~\ref{Beta} introduces the selection criteria of $\beta$, and investigates the effect of $\beta$ on the stability and accuracy.

\subsection{Regularization Parameter $\alpha$} \label{RegularAlpha}

\begin{figure}
	\centering
		\includegraphics[scale=.55]{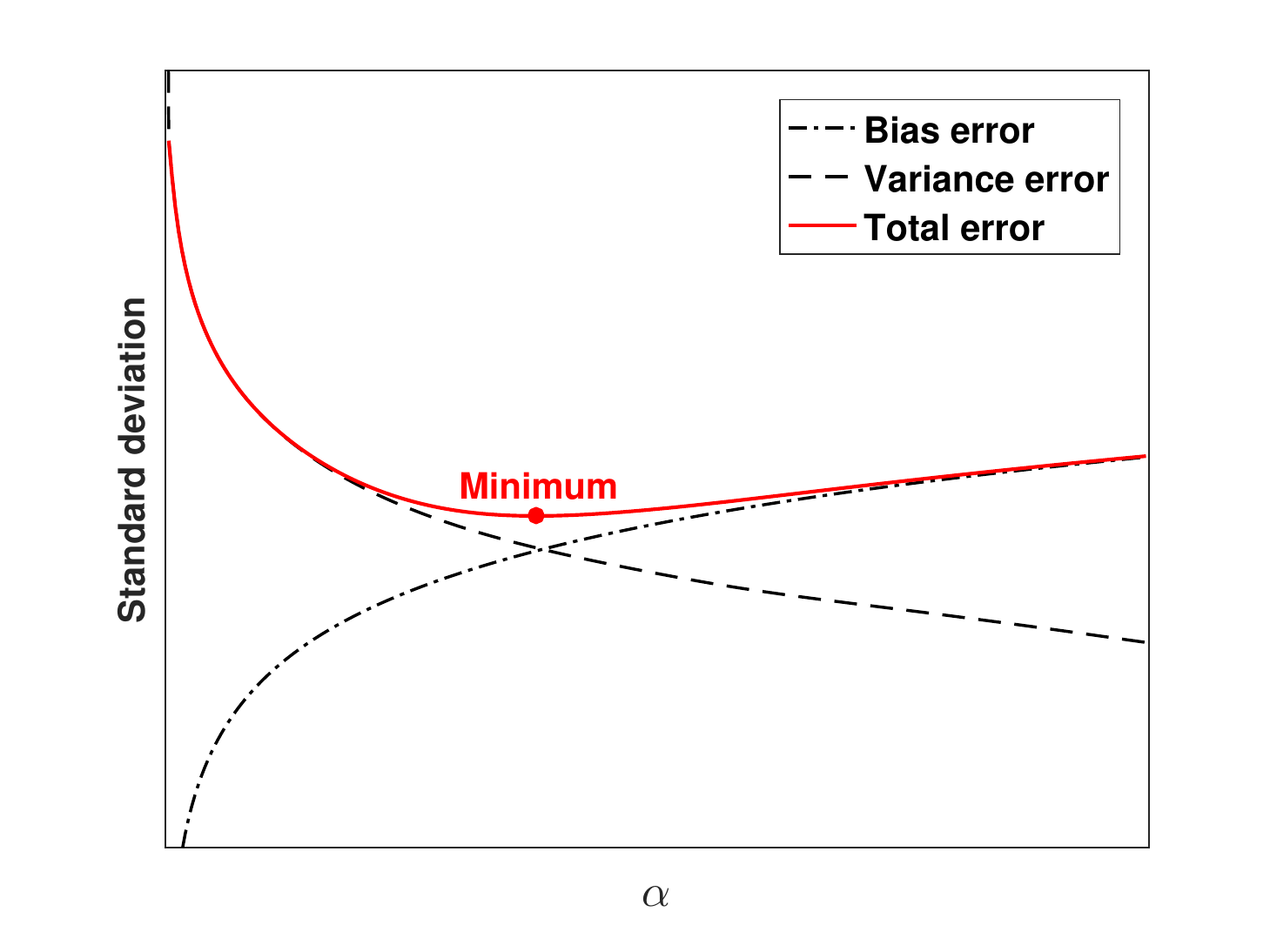}
	\caption{Bias-variance trade-off conceptual plot in the regularization method.}
	\label{ErrorType}
\end{figure}

The optimal regularization parameter $\alpha$ is unknown and can be mathematically obtained by,
     \begin{equation} \label{Terror}
     \mathscr{T}^2(\mb{q}_{\alpha},\sigma) = \frac{1}{M}\sum_{m=1}^{M}\left[ \tilde{\mb{q}}^m - \mb{q}_{\alpha}^m(\sigma) \right]^T\left[ \tilde{\mb{q}}^m - \mb{q}_{\alpha}^m(\sigma) \right],
     \end{equation}
where $\mathscr{T}$ is the total error, and $\tilde{\mb{q}}^m$ is a true heat flux component vector. 
However, using the residual sum of squares in Eq.~\ref{Terror} is impractical in real experimental environments because the true heat flux vector $\tilde{\mb{q}}^m$ is unknown.
Therefore, the optimal $\alpha$ is estimated using the discrepancy principle in Eq.~\ref{Morozov} in general.
Monte-Carlo methods~\cite{McDonald1975, Honerkamp1990}, which are well known computational algorithms to obtain numerical results, are then used to specify the parameter $\alpha$ in Eqs~\ref{Morozov} and~\ref{Terror}.
However, the high computational demands of the Monte-Carlo methods including repeated parametric sampling may not be suitable for the proposed iterative scheme to select the hybrid parameters.
Note that the discrepancy principle in Eq~\ref{Morozov} is a temperature-based ridge estimator.
It may yield lower accuracy for heat flux estimation than for temperature estimation.
These two limitations could be addressed by decomposing the total error in Eq.~\ref{Terror}, as follows~\cite{Beck1985}. 
     \begin{equation} \label{TDV}
     \mathscr{T}^2(\mb{q}_{\alpha},\sigma) = \mathscr{D}^2(\Delta\mb{q}, \alpha) + \mathscr{V}^2(\sigma, \alpha),
     \end{equation}
where $\mathscr{D}$ is a bias (or deterministic) error and is a function of $\alpha$ and the heat flux change, and
$\mathscr{V}$ is a variance (or stochastic) error, which is a function of $\alpha$ and the measurement noise.
This implies that the total error could be decomposed by the bias and variance errors.
Figure~\ref{ErrorType} conceptually shows a plot of $\mathscr{T}$, $\mathscr{V}$, and $\mathscr{D}$ as $\alpha$ is varied.
By separating the total error in Eq.~\ref{TDV}, the input data (measurement noise, $\sigma$) and output data (heat flux solution, $\mb{q}_{\alpha}$) can be treated independently.
Making it possible to improve the accuracy and efficiency of the conventional ridge estimator.
In the following subsections, we discuss how to calculate each error efficiently.

\subsubsection{Bias error}   \label{Bias}

The bias error of the heat flux is defined as 
     \begin{equation} \label{Derror_1}
     \mathscr{D}^2 = \frac{1}{M}\sum_{m=1}^{M}\left(\delta \mb{q}_{bias}^{m}\right)^T\left(\delta \mb{q}_{bias}^{m}\right),
     \end{equation}
in which $\delta \mb{q}_{bias}^{m}$ is the bias error vector of the heat flux and is a function of $\alpha$ and heat flux change.
To calculate the bias error vector, the gradient of Eq.~\ref{Solver} is set to zero, which yields
     \begin{equation} \label{GradSreg}
     \alpha\Delta\mb{q}_{\alpha}^{m+1} = \mb{X}_s^T \left[\mb{Y}^{m+1} - \mb{L}_s\mb{T}_{\alpha}^{m+1}\right],  
     \end{equation}
where $\mb{L}_s\mb{T}_{\alpha}^{m+1}$ can be defined as
     \begin{equation} \label{Tb}
     \mb{L}_s\mb{T}_{\alpha}^{m+1} = \mb{Y}^{m+1} - \delta\mb{T}_{bias}^{m+1}.
     \end{equation}

Here $\delta\mb{T}_{bias}^{m+1}$ is th bias error vector of the temperature corresponding to the sensor positions and is a function of $\alpha$.
When $\alpha$ is zero, $\delta\mb{T}_{bias}^{m+1}$ becomes a zero vector, which indicates that $\mb{L}_s\mb{T}_{\alpha}^{m+1}$ is over-fitted to $\mb{Y}^{m+1}$.
Substituting Eq.~\ref{Tb} into Eq.~\ref{GradSreg} gives
     \begin{equation} \label{Tbias}
     \delta\mb{T}_{bias}^{m+1} = \alpha \left(\mb{X}_s\mb{X}_s^T\right)^{+}\mb{X}_s \Delta\mb{q}_{\alpha}^{m+1},  
     \end{equation}
where $\left(\mb{X}_s\mb{X}_s^T\right)^{+}$ is the Moore-Penrose generalized inverse of $\mb{X}_s\mb{X}_s^T$.
Equation~\ref{Tbias} shows that the bias error of the temperature is a function of $\alpha$ and $\Delta\mb{q}_{\alpha}^{m+1}$.

\begin{figure}
\centering        
        \includegraphics[scale=.3]{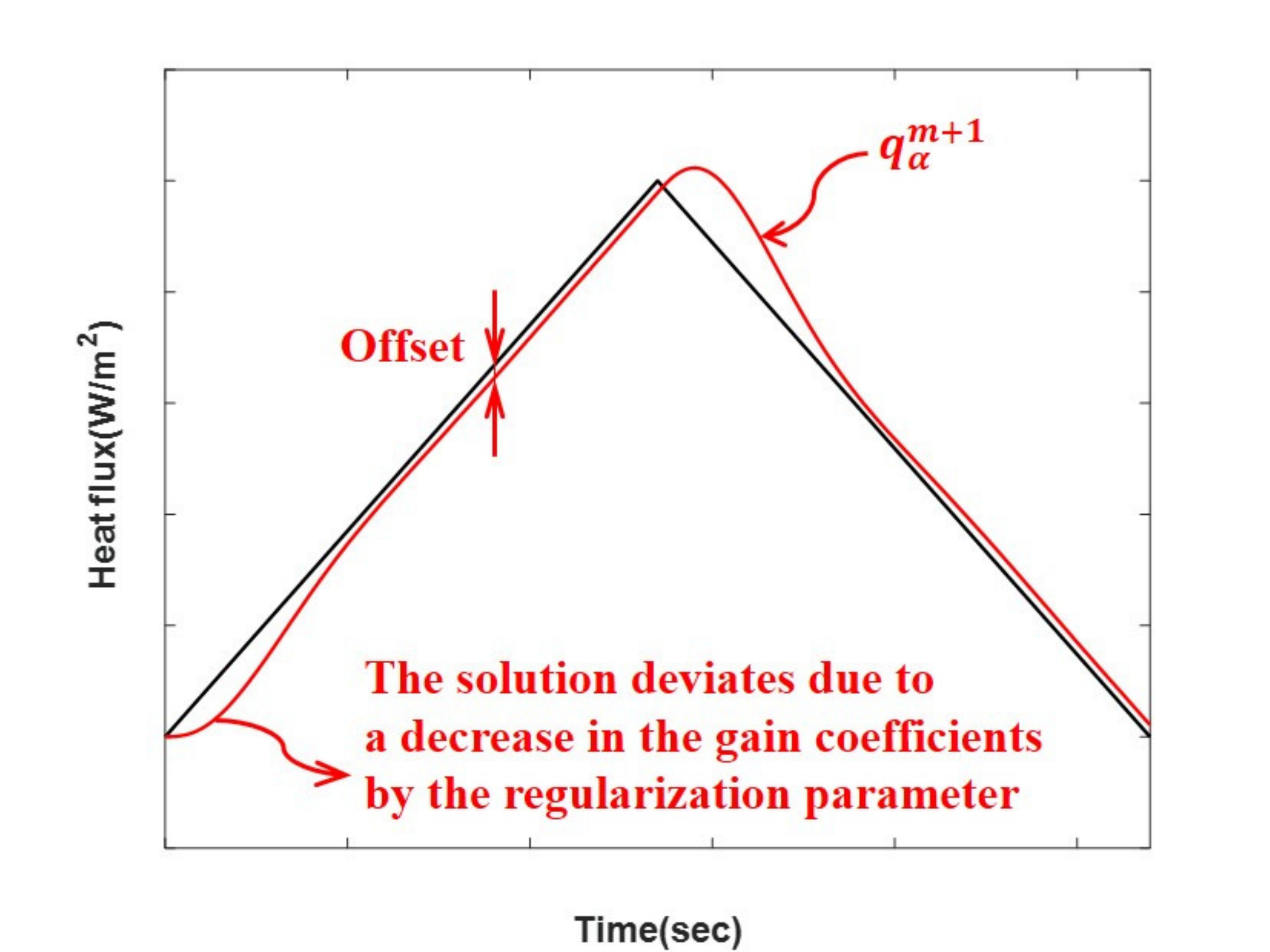}
\caption{Illustration of the bias error caused by excessive regularization.}
    \label{BiasError}
    \end{figure}

The heat flux solution tends to deviate initially and to become parallel to the exact solution (Figure~\ref{BiasError})---
these are called the deviation error and offset error, respectively.
The deviation error is then neglected when the sensors are attached to certain positions with high sensitivity coefficients.
A sensitivity analysis is essential for this purpose, to find adequate sensor locations in general~\cite{Beck1985}.
In this work, we assume the deviation error is sufficiently small; then, $\Delta\mb{q}_{\alpha}^{m+1}$ can be defined by considering only the offset error, as follows.
     \begin{equation} \label{qoffset}
     \Delta\mb{q}_{\alpha}^{m+1}=\Delta\mb{q}_{exact}^{m+1},  
     \end{equation}
where $\Delta\mb{q}_{exact}^{m+1}$ is the exact heat flux change and is independent of $\alpha$.
This implies that it is unnecessary to calculate the heat flux change for all discretized $\alpha$ to obtain $\delta\mb{T}_{bias}^{m+1}$.

In addition, the maximum bias error could be considered to avoid the computational requirement of  $\Delta\mb{q}_{exact}^{m+1}$ at each time--step; then, $\Delta\mb{q}_{exact}^{m+1}$ of Eq.~\ref{qoffset} can be defined as,
     \begin{equation} \label{qoffset2}
     \Delta\mb{q}_{exact}^{m+1} \coloneqq \Delta\mb{q}_{max} \qquad \text{for} \: \left\lVert \Delta\mb{q}_{max} \right\lVert= \max\limits_{1\le m \le M} \left\lVert \Delta\mb{q}_{exact}^{m} \right\lVert,
     \end{equation}
where $\left\lVert\cdot\right\lVert$ denotes the $L^2$-norm (or Euclidean norm).
Substituting Eqs.~\ref{qoffset} and~\ref{qoffset2} into Eq.~\ref{Tbias} yields
     \begin{equation} \label{TbiasNew}
      \delta\mb{T}_{bias}^{m+1} \coloneqq \delta\mb{T}_{\textit{offset}} = \alpha \left(\mb{X}_s\mb{X}_s^T\right)^{+}\mb{X}_s \Delta\mb{q}_{max}.
     \end{equation}

Equation~\ref{TbiasNew} shows that $\alpha$ is only a variable; thus, the computational cost of $\delta\mb{T}_{\textit{offset}}$ is much lower than that of $\delta\mb{T}_{bias}^{m+1}$ in Eq.~\ref{Tbias}.
If $\Delta\mb{q}_{max}$ is difficult to define, the mean value of $\Delta\mb{q}_{\alpha}^{m+1}$ may be used to simplify $\Delta\mb{q}_{exact}^{m+1}$. 

Considering the bias error ($ \delta\mb{T}_{bias}^{m+1}$) in the discrete system of Eq.~\ref{FEM}, we have
     \begin{equation} \label{qoffset3}
     \left(\sum_{j=1}^N \delta q_{bias}^{(j)^{m+1}}\bar{\mb{q}}^{(j)}\right) = \mb{C} \delta\mb{\dot{T}}_{bias}^{m+1} + \mb{K}\delta\mb{T}_{bias}^{m+1},  
     \qquad \frac{\partial \mb{q}_g}{\partial \mb{q}} =  \left[ \bar{\mb{q}}^{(1)}, \cdots, \bar{\mb{q}}^{(N)} \right],
     \end{equation} 
where $\delta q_{bias}^{(j)^{m+1}}$ is a bias error of the heat flux component specified on the $j$-th boundary.
$\left(\partial \mb{q}_g/ \partial \mb{q}\right)$ is defined in Eq.~\ref{UnitQ}. Here, $\bar{\mb{q}}^{(j)}$ is a unit heat flux vector specified on the $j$-th boundary.
Considering only the offset error, Eq.~\ref{qoffset3} can be rewritten as
     \begin{equation} \label{qoffset4}
     \left(\sum_{j=1}^N \delta q_{\textit{offset}}^{(j)}\bar{\mb{q}}^{(j)}\right) = \left(\lambda\mb{C} + \mb{K}\right)\delta\mb{T}_{\textit{offset}} \:, 
     \end{equation} 
where $\lambda$ is a thermal eigenvalue, which is unknown and generally negative in thermal problems.
To handle the unknown $\lambda$, we employ O’Callahan’s approach~\cite{OCallahan}, and the following relation with no loading condition is obtained:
     \begin{equation} \label{LambdaT}
     \lambda\delta\mb{T}_{\textit{offset}} \approx \mb{C}^{-1}\mb{K}\delta\mb{T}_{\textit{offset}} \:.
     \end{equation}

The dimensions of $\mb{K}$ and $\delta\mb{T}_{\textit{offset}}$ are $\left(N_g \times N_g\right)$ and $\left(I \times 1\right)$, respectively.
The number of sensors is much smaller than the number of total DOFs ($I \ll N_g$); therefore, most DOFs are not measured. 
Equation~\ref{qoffset4} with the assumption in Eq.~\ref{LambdaT} is then rewritten as,
     \begin{equation} \label{qoffset5}
     \left(\sum_{j=1}^N \delta q_{\textit{offset}}^{(j)}\begin{bmatrix} \bar{\mb{q}}_s^{(j)} \\ \bar{\mb{q}}_v^{(j)} \end{bmatrix} \right)=
     2\begin{bmatrix} \mb{K}_{ss} & \mb{K}_{sv} \\ \mb{K}_{vs} & \mb{K}_{vv} \end{bmatrix} \begin{bmatrix} \delta\mb{T}_{\textit{offset}} \\ \delta\mb{T}_{v} \end{bmatrix},
     \qquad \mb{K} = \begin{bmatrix} \mb{K}_{ss} & \mb{K}_{sv} \\ \mb{K}_{vs} & \mb{K}_{vv} \end{bmatrix},
     \qquad \bar{\mb{q}}^{(j)} = \begin{bmatrix} \bar{\mb{q}}_s^{(j)} \\ \bar{\mb{q}}_v^{(j)} \end{bmatrix}, 
     \end{equation}
where the subscript $v$ denotes virtual parts, which are quantities with respect to the unmeasured DOFs.
$\delta\mb{T}_{v}$ is a virtual temperature error vector of the unmeasured positions.
The two resulting equations in Eq.~\ref{qoffset5} yield 
     \begin{subequations} \label{Condensed}
     \begin{gather}
     \label{qoffset6}
     \left(\sum_{j=1}^N \delta q_{\textit{offset}}^{(j)}\hat{\mb{q}}_s^{(j)}\right) = \hat{\mb{K}}\delta\mb{T}_{\textit{offset}},\\
     \label{qoff_comps}
     \hat{\mb{q}}_s^{(j)} = \bar{\mb{q}}_s^{(j)} + \mb{\Psi}_{vs}^T\bar{\mb{q}}_v^{(j)}, \qquad \hat{\mb{K}} = 2\left(\mb{K}_{ss} + \mb{K}_{sv}\mb{\Psi}_{vs}\right), \qquad \mb{\Psi}_{vs} = -\mb{K}_{vv}^{-1}\mb{K}_{vs},
     \end{gather}
     \end{subequations}

Then, multiplying by $\delta\mb{T}_{\textit{offset}}^T$ on both sides yields
     \begin{subequations} \label{CondensedQ}
     \begin{gather}
     \label{qoffset7}
     \delta \mb{q}_{\textit{offset}}^T \mb{q}_{\delta T} = \delta\mb{T}_{\textit{offset}}^T\hat{\mb{K}}\delta\mb{T}_{\textit{offset}},\\
     \label{qoff_compsQoff}
     \delta \mb{q}_{\textit{offset}} = \left[ \delta q_{\textit{offset}}^{(1)}, \cdots, \delta q_{\textit{offset}}^{(N)} \right]^T, 
     \qquad \mb{q}_{\delta T} = \left[ \delta\mb{T}_{\textit{offset}}^T\hat{\mb{q}}_s^{(1)} , \cdots, \delta\mb{T}_{\textit{offset}}^T\hat{\mb{q}}_s^{(N)}\right]^T.
     \end{gather}
     \end{subequations}

Equation~\ref{qoffset7} can be restated by the Cauchy–Schwarz inequality~\cite{Paulsen2002}, which gives
     \begin{equation} \label{qoffset9}
     \left\| \delta \mb{q}_{\textit{offset}}\right\| \Big \|\mb{q}_{\delta T} \Big \| \geq \left\| \delta\mb{T}_{\textit{offset}}^T\hat{\mb{K}}\delta\mb{T}_{\textit{offset}}\right\|. 
     \end{equation}

Using Eq.~\ref{qoffset9}, the bias error of the heat flux in Eq.~\ref{Derror_1} can be redefined as,
     \begin{equation} \label{DerrorNew}
     \mathscr{D}^2 = \left\| \delta \mb{q}_{\textit{offset}}\right\|^2 \geq \left\| \delta\mb{T}_{\textit{offset}}^T\hat{\mb{K}}\delta\mb{T}_{\textit{offset}}\right\|^2 \bigm/ \Big \|\mb{q}_{\delta T} \Big \|^2,
     \qquad \delta\mb{T}_{\textit{offset}} = \alpha \left(\mb{X}_s\mb{X}_s^T\right)^{+}\mb{X}_s \Delta\mb{q}_{max},
     \end{equation}
where $\delta \mb{q}_{\textit{offset}}$ is a time-independent function unlike $\delta \mb{q}_{bias}^{m}$ in Eq.~\ref{Derror_1}.
Therefore, the computational efficiency is much better for Eq.~\ref{DerrorNew} than for Eq.~\ref{Derror_1}.
The deviation error is neglected here to obtain better computational efficiency; however, it may have a limitation in certain cases with step or constant heat fluxes.
This is because the offset error is zero when there is no heat flux change.

\subsubsection{Variance error}   \label{Variance}
The variance error is mainly caused by the propagation and accumulation of the measurement error. The variance error of the heat flux is defined as 
     \begin{equation} \label{Verror}
     \mathscr{V}^2 = \frac{1}{M}\sum_{m=1}^{M}\left(\delta \mb{q}_{var}^{m}\right)^T\left(\delta \mb{q}_{var}^{m}\right),
     \end{equation}
where $\delta \mb{q}_{var}^{m}$ is a variance error vector of the heat flux, and is a function of $\alpha$ and measurement noise.
To calculate the propagation and accumulation of the measurement error over time, Eq.~\ref{MNEq_r} is written as,
     \begin{equation} \label{MNESum_error}
     \Delta \mb{q}_{\alpha}^{m+1} = \mb{G}_{\alpha} \left[\mb{Y}^{m+1} - 
     \mb{L}_s\left\{\mb{A}^{m+1}\mb{T}^0 + \left(\sum_{p=0}^{m}\mb{A}^p {U}\mb{q}_g^0 \right) + \left(\sum_{i=0}^{m} \sum_{k=0}^{i+1}\mb{A}^k \mb{X}\Delta \mb{q}_{\alpha}^{m-i}\right)\right\} \right].
     \end{equation}

If the measurement error follows a normal distribution, $\mb{Y}^{m+1}$ can be defined as,
     \begin{equation} \label{Y}
     \mb{Y}^{m+1} = \tilde{\mb{Y}}^{m+1} + \delta\mb{Y}^{m+1}, \qquad \delta\mb{Y}^{m+1} \sim \mathcal{N}(\mb{0},\boldsymbol{\sigma}^2) \quad \text{for $m=1,\cdots,M-1$,}
     \end{equation}
where $\tilde{\mb{Y}}^{m+1}$ and $\delta\mb{Y}^{m+1}$ are the exact temperature and measurement error vectors, respectively.
$\delta\mb{Y}^{m+1}$ is normally distributed with zero-mean and standard deviation, denoted by $\delta\mb{Y}^{m+1} \sim \mathcal{N}(\mb{0},\boldsymbol{\sigma}^2)$,
and $\boldsymbol{\sigma}$ is a $I$-dimensional standard deviation (i.e., measurement noise level of the sensors) vector.

In general, the measurement noise level of temperature sensors is much larger than the numerical error.
The initial boundary condition error is negligible if the thermal boundary condition is clearly defined with a well-established experimental environment. 
Then, Eq.~\ref{MNESum_error} could be defined using only the measurement error:
     \begin{equation} \label{MNESum_error2}
     \delta\Delta \mb{q}_{var}^{m+1} 
     = \mb{G}_{\alpha} \left[\delta\mb{Y}^{m+1} - \mb{L}_s\left( \sum_{i=0}^{m} \left(\mb{I} - \mb{A}^{i+2}\right)\left(\mb{I} - \mb{A}\right)^{-1} \mb{X} \delta\Delta \mb{q}_{\alpha}^{m-i}\right) \right],\\ 
     \end{equation}

In a similar manner, Eq.~\ref{MNEq_r} could also be written as,
     \begin{equation} \label{MNESum_error2_17anew}
     \delta\Delta \mb{q}_{\alpha}^{m-i} = \mb{G}_{\alpha} \delta\mb{Y}^{m-i}.
     \end{equation}

Substituting Eq.~\ref{MNESum_error2_17anew} into Eq.~\ref{MNESum_error2}, we have
     \begin{equation} \label{MNESum_error2_new}
     \delta\Delta \mb{q}_{var}^{m+1} 
     = \mb{G}_{\alpha} \left[\delta\mb{Y}^{m+1} - \mb{L}_s\left( \sum_{i=0}^{m} \left(\mb{I} - \mb{A}^{i+2}\right)\left(\mb{I} - \mb{A}\right)^{-1} \mb{X} \mb{G}_{\alpha} \delta\mb{Y}^{m-i}\right) \right],
     \end{equation}
where $\delta\Delta \mb{q}_{var}^{m+1}$ is a time-dependent variance error of heat flux change. 
Equation~\ref{MNESum_error2_new} shows that $\delta\Delta \mb{q}_{var}^{m+1}$ is a function of the measurement error and $\alpha$.

The summation term in Eq.~\ref{MNESum_error2_new} represents the propagation and accumulation of an error, and it entails high computational resources.
To address this problem, we assume $\delta\mb{Y}^{m+1}$ as a regular oscillation:
     \begin{equation} \label{SigmaU}
     \delta\mb{Y}^{m+1} = \left(-1\right)^{(m+1)}\boldsymbol{\sigma},
     \end{equation}
where the superscript $( \cdot )$ is a power (or exponent).
By substituting Eq.~\ref{SigmaU} into Eq.~\ref{MNESum_error2_new}, we can obtain a time-independent variance error of the heat flux change:
     \begin{subequations} \label{MNESum_error3}
     \begin{gather} 
     \label{}
     \delta\Delta \mb{q}_{var} = \mb{G}_{\alpha} \left[\mb{I} + \mb{L}_s \left[\left\{\mb{I} + \mb{A}^{(M_a+1)}\right\} 
     \left(\mb{I} + \mb{A}\right)^{-1} \mb{A}^{(2)} - \mb{I}\right]\left(\mb{I} - \mb{A}\right)^{-1} \mb{X} \mb{G}_{\alpha} \right]\boldsymbol{\sigma}, \qquad(\text{when $M_a$ is even}), \\
     \label{}
     \delta\Delta \mb{q}_{var} = \mb{G}_{\alpha} \left[\mb{I} + \mb{L}_s \left\{\mb{I} - \mb{A}^{(M_a+1)}\right\} 
     \left(\mb{I} + \mb{A}\right)^{-1} \mb{A}^{(2)}\left(\mb{I} - \mb{A}\right)^{-1} \mb{X} \mb{G}_{\alpha} \right]\boldsymbol{\sigma}, \:\qquad\qquad (\text{when $M_a$ is odd}),
     \end{gather}
     \end{subequations}
where $M_a$ is a time step at the moment when the norm of $\mb{A}^{(M_a+1)}$ is less than one.
Using Eq.~\ref{MNESum_error3}, Eq.~\ref{Verror} can be redefined as,
     \begin{equation} \label{VerrorNew}
     \mathscr{V}^2 = \left(\delta\mb{q}_{var}\right)^T\left(\delta\mb{q}_{var}\right), \qquad \delta\mb{q}_{var} = \xi(\alpha) \delta\Delta \mb{q}_{var},
     \end{equation}
where $\delta\mb{q}_{var}$ is a time-independent variance error of the heat flux, and $\xi(\alpha)$ is the ratio of $\delta\mb{q}_{var}$ to $\delta\Delta \mb{q}_{var}$.
If the measurement noise is defined as Eq.~\ref{SigmaU}, $\xi(0)$ is exactly 0.5.
As $\alpha$ increases, $\xi(\alpha)$ increases slightly.
If the measurement noise follows the normal distribution, $\xi(0)$ is close to, and not exactly, 0.5.
However, as $\alpha$ increases, $\xi(\alpha)$ rapidly increases because of the uncertainty of the measurement.
To compensate the assumption in Eq.~\ref{SigmaU}. $\xi(\alpha)$ is defined using the measurement uncertainty, as follows.
     \begin{equation} \label{xi}
     \xi(\alpha) = \sqrt{\frac{\sum_{m=1}^{M_{s}}|\bar{\mb{q}}_{\alpha}^{m}|^2}{\sum_{m=1}^{M_{s}}|\Delta \bar{\mb{q}}_{\alpha}^{m}|^2}},
     \end{equation}
where $\bar{\mb{q}}_{\alpha}^m$ and $\Delta \bar{\mb{q}}_{\alpha}^m$ are the heat flux and heat flux change vectors, respectively, calculated by $\delta \mb{Y}^m$, which is the measurement noise data with zero--mean.
If $\sigma$ is known, $\delta\mb{Y}^m$ can be easily generated. 
$t_{M_s}$ is much smaller than the total measurement time $t_M$. 
This is because the expected mean value of $\delta\mb{Y}^m$ is zero by the standard statistical assumptions, and thus $\xi(\alpha)$ can be obtained in a short-term trend.

Finally, the total error of the heat flux is calculated as the sum of $\mathscr{D}^2$  and $\mathscr{V}^2$ defined in Eqs.~\ref{DerrorNew} and~\ref{VerrorNew}, respectively.
 It is important that Eqs.~\ref{DerrorNew} and~\ref{VerrorNew} do not depend on the heat flux history. 
Therefore, it is possible to estimate the optimal $\alpha$ more efficiently than the conventional methods.

\subsection{Parameter $\beta$}  \label{Beta}

In the previous section, the optimal selection process of $\alpha$ was described.  
However, $\alpha$ is a penalty value, indicating that the bias error increases when $\alpha$ increases. 
This may cause excessive fluctuation of the inverse solution at the initial time--steps and/or a certain time step with abrupt heat flux change.
This could be addressed by controlling $\beta$ with less penalty than $\alpha$.
The parameter $\beta$ should be selected such that the magnitude of the gain coefficient matrix is minimal, in which case the stability conditions should also be satisfied.
Therefore, in this section, we formulate the stability condition for $\beta$.

First, stability can be evaluated in terms of the inverse solution in Eq.~\ref{Test}:
     \begin{equation} \label{T_est}
     \mb{T}_{\alpha}^{m+1} = \bar{\mb{T}}_{\alpha}^{m+1} + \mb{X}\Delta \mb{q}_{\alpha}^{m+1}. \\
     \end{equation}

Substituting Eqs.~\ref{TaylorSetT} and~\ref{MNEq_r} into Eq.~\ref{T_est} gives
     \begin{equation} \label{Tstability}
     \mb{T}_{\alpha}^{m+1} = \mb{E}_{\alpha}\mb{A}\mb{T}_{\alpha}^m + \mb{E}_{\alpha}\mb{U}\mb{q}_{g,\alpha}^m + \mb{X}\mb{G}_{\alpha}\mb{Y}^{m+1}, \qquad \mb{E}_{\alpha} = \mb{I} - \mb{X}\mb{G}_{\alpha}\mb{L}_s, \\
     \end{equation}
where $\mb{q}_{g,\alpha}^m$ is a global nodal heat flux vector determined by the inverse process with the regularization method.
Equation~\ref{Tstability} is clearly decoupled into two-parts: the inverse solution (temperature and heat flux) at the current time--step and measurement data.
$\mb{E}_{\alpha}\mb{A}\mb{T}_{\alpha}^m$ and $\mb{E}_{\alpha}\mb{U}\mb{q}_{g,\alpha}^m$ are the inverse solutions at the current time--step, and $\mb{X}\mb{G}_{\alpha}\mb{Y}^{m+1}$ denotes the measurement data.
Stability depends on the initial error of the solution with respect to $\mb{E}_{\alpha}\mb{A}\mb{T}_{\alpha}^m$ and $\mb{E}_{\alpha}\mb{U}\mb{q}_{g,\alpha}^m$
, and the input data of $\mb{X}\mb{G}_{\alpha}\mb{Y}^{m+1}$ (i.e., measurement noise data) is also a significant factor influencing the stability of the inverse solution.
Therefore, both the initial and the input data errors should be well--controlled to realize stable inverse solutions.

First, for the initial error bound, the error of $\mb{T}_{\alpha}^{m+1}$ satisfies the following equations
     \begin{subequations} \label{InitialError_new}
     \begin{gather}
     \delta \mb{T}_{\alpha}^{m+1} \coloneqq \left(\mb{E}_{\alpha}\mb{A}\right)^{(m+1)}\delta\mb{T}_{\alpha}^0, \qquad  \left(\text{when}  \left\|\mb{E}_{\alpha}\mb{A}\right\| > \left\|\mb{E}_{\alpha}\mb{U}\right\| \right),\\
     \delta \mb{T}_{\alpha}^{m+1} \coloneqq \left(\mb{E}_{\alpha}\mb{U}\right)^{(m+1)}\delta\mb{q}_{g,\alpha}^0, \quad\:\:\:  \left(\text{when}  \left\|\mb{E}_{\alpha}\mb{A}\right\| < \left\|\mb{E}_{\alpha}\mb{U}\right\| \right),                     
     \end{gather}
     \end{subequations}
where $\mb{E}_{\alpha}\mb{A}$ and $\mb{E}_{\alpha}\mb{U}$ are amplification matrices. 
The stability of the discretization methods requires certain conditions to ensure the uniform boundedness of the arbitrary powers of the amplification matrix.
Therefore, the spectral radius of $\mb{E}_{\alpha}\mb{A}$ and $\mb{E}_{\alpha}\mb{U}$ must be smaller than 1 to bound the error, which is the well-known von--Neumann condition~\cite{Hirsch1988}.
By definition, the spectral radius is the largest absolute value of its eigenvalues, which can be written as
     \begin{equation} \label{SpecR}
     \rho\left(\mb{E}_{\alpha}\mb{A}\right) \coloneqq \text{max}\left\{|\bar{\lambda}_1|,\cdots,|\bar{\lambda}_{N_g}|\right\} \leq 1,
     \end{equation}
where $\bar{\lambda}$ is an eigenvalue of $\mb{E}_{\alpha}\mb{A}$. 

The inverse solution does not divergence regardless of $\mb{X}\mb{G}_{\alpha}\mb{Y}^{m+1}$, if the stability condition for $\beta$ given in Eq.~\ref{SpecR} is satisfied.
However, if the stability condition for $\mb{X}\mb{G}_{\alpha}\mb{Y}^{m+1}$ is not sufficient, the amplitude of the temperature oscillation is much greater than the measurement errors.
Consequently, it is difficult to distinguish if the changes in the boundary heat flux or input errors are the main reason for the fluctuation of the estimated temperature.
Therefore, we should revisit Section~\ref{RegularAlpha} to select $\alpha$  again to handle the instability of $\mb{X}\mb{G}_{\alpha}\mb{Y}^{m+1}$.
This is the main motivation to propose an iterative algorithm to select $\alpha$ and $\beta$ in this work.
The details are described in the following section.

\subsection{Algorithm to Select $\alpha$ and $\beta$}  \label{Algorithm0}

In this section, the selection criteria for $\alpha$ and $\beta$ are presented.
The regularization method provides improved accuracy and stability in exchange for a tolerable amount of bias (see the bias--variance trade-off in Figure~\ref{ErrorType}).
At this time, if $\beta$ is properly considered, the variance error can be reduced with less penalty than $\alpha$.
Therefore, selecting $\alpha$ properly with $\beta$ reduces the amount of bias error that needs to be exchanged, thereby reducing the total error.
The specific criteria are discussed in Sections~\ref{RegularAlpha} and~\ref{Beta}, and the overall process of selecting hybrid parameters is described in Algorithm~\ref{Algorithm}. 
In Algorithm~\ref{Algorithm}, $N_{iter}$ denotes the number of iterations.
Conventional methods require $\sigma$, $\Delta \alpha$, $\Delta \beta$, $\Delta t$, and $\mb{Y}$ as input data. 
On the other hand, in the proposed algorithm, $\Delta\mb{q}_{max}$ and $t_{M_s}$ replace $\mb{Y}$.
Note that $\xi(\alpha)$ is obtained in the form of oscillation, because $\xi(\alpha)$ is calculated with the noise data.
To handle this, a fitted curve can be used.
Note further that the algorithm converges fast with $N_{iter} = 2$ in the numerical problems considered in this work.
\begin{algorithm}
\caption{}
\label{Algorithm}
\renewcommand{\algorithmicrequire}{\textbf{Input:}}
\renewcommand{\algorithmicensure}{\textbf{Output:}}
\def\NoNumber#1{{\def\alglinenumber##1{}\State #1}\addtocounter{ALG@line}{-1}}
\begin{algorithmic}[1]
   \Require{$\mb{C}$, $\mb{K}$, $\mb{X}$, $\Delta \mb{q}_{max}$, $\sigma$, $\Delta \alpha$, $\Delta \beta$, $\Delta t$, $t_{M_s}$}
   \Ensure{$\alpha$, $\beta$ }
   \State  $\alpha=0$, \quad $\beta=1$,  \quad $\xi(\alpha)$                                             \Comment{Initial setting (see Eq.~\ref{xi})} 
   \For{$k = 1$ to $N_{iter}$}
       \State  $\mathscr{T}^2 = \mathscr{T}_{old}^{2}= \left\|\delta\mb{T}_{\textit{offset}}^T\hat{\mb{K}}\delta\mb{T}_{\textit{offset}}\right\|^2 \bigm/ \Big \|\mb{q}_{\delta T} \Big \|^2 + \xi(\alpha)^2\left(\delta\Delta\mb{q}_{var}\right)^T\left(\delta\Delta\mb{q}_{var}\right)$                                                          \Comment{(see Eqs.~\ref{DerrorNew} and~\ref{VerrorNew})} 
       \While{$ \mathscr{T}_{old}^{2}  \ge \mathscr{T}^2$}
            \State $\alpha = \alpha + \Delta \alpha$
                                                       
            \State $\mathscr{T}_{old}^{2}=\mathscr{T}^2$
            \State $\mathscr{T}^2 = \left\|\delta\mb{T}_{\textit{offset}}^T\hat{\mb{K}}\delta\mb{T}_{\textit{offset}}\right\|^2 \bigm/ \Big \|\mb{q}_{\delta T} \Big \|^2 + \xi(\alpha)^2\left(\delta\Delta\mb{q}_{var}\right)^T\left(\delta\Delta\mb{q}_{var}\right)$
       \EndWhile
       \State $\alpha = \alpha - \Delta \alpha$
       \If{$k==N_{iter}$}
            \State break   \Comment{Output: $\alpha$ and $\beta$} 
       \EndIf
       \State $\beta_0=\beta$                                                                                                  \Comment{For stopping criteria of $\beta$}  
       \State $\mb{G}_{\alpha} = \mb{G}_{\alpha}^{old}=\left(\mb{X}_s^T\mb{X}_s + \alpha\mb{I}\right)^{-1}\mb{X}_s^T$          \Comment{(see Eq.~\ref{MNEq_r})}   
       \While{$||\mb{G}_{\alpha}^{old}||\ge||\mb{G}_{\alpha}||$ and $\rho(\mb{E}_{\alpha}\mb{A})<1$}
            \State $\beta=\beta - \Delta\beta$  
            \State $\mb{G}_{\alpha}^{old} = \mb{G}_{\alpha}$
            \State $\mb{G}_{\alpha} =\left(\mb{X}_s^T\mb{X}_s + \alpha\mb{I}\right)^{-1}\mb{X}_s^T$
            \State $\mb{E}_{\alpha}\mb{A} = (\mb{I} - \mb{X}\mb{G}_{\alpha}\mb{L}_s)\mb{A}$              \Comment{(see Eqs.~\ref{TgA} and~\ref{Tstability})} 
       \EndWhile
       \State $\beta = \beta + \Delta\beta$  
       \If{$\beta == 1$ or $\beta = \beta_0$}
            \State break   \Comment{Output: $\alpha$ and $\beta$} 
       \EndIf    
       
       \State $\mb{A}=\mb{A}(\beta)$,\quad $\mb{U}=\mb{U}(\beta)$, \quad $\mb{X}=\mb{X}(\beta)$          \Comment{Update the component matrices} 
   \EndFor  
   \State \textbf{end}
\end{algorithmic}
\end{algorithm}

\section{Numerical Examples} \label{NumericalExamples}
This Section evaluates the performance of the algorithm with two numerical models: 1D bar and 2D plate.
In Section~\ref{Verify_a}, the proposed estimator is evaluated with changes in the thermal diffusivity and the heat flux profile.
In Section~\ref{Verify_b}, we verify the stability conditions for $\beta$ and demonstrate that the bias error decreases when $\beta$ is considered.
In Sections~\ref{perform1D} and~\ref{EX2D}, the hybrid parameter selection algorithm is validated using the 1D and 2D examples,respectively, with various conditions.

\subsection{One-dimensional problem}  \label{EX1D}
\begin{figure}
\centering        
        \includegraphics[scale=.6]{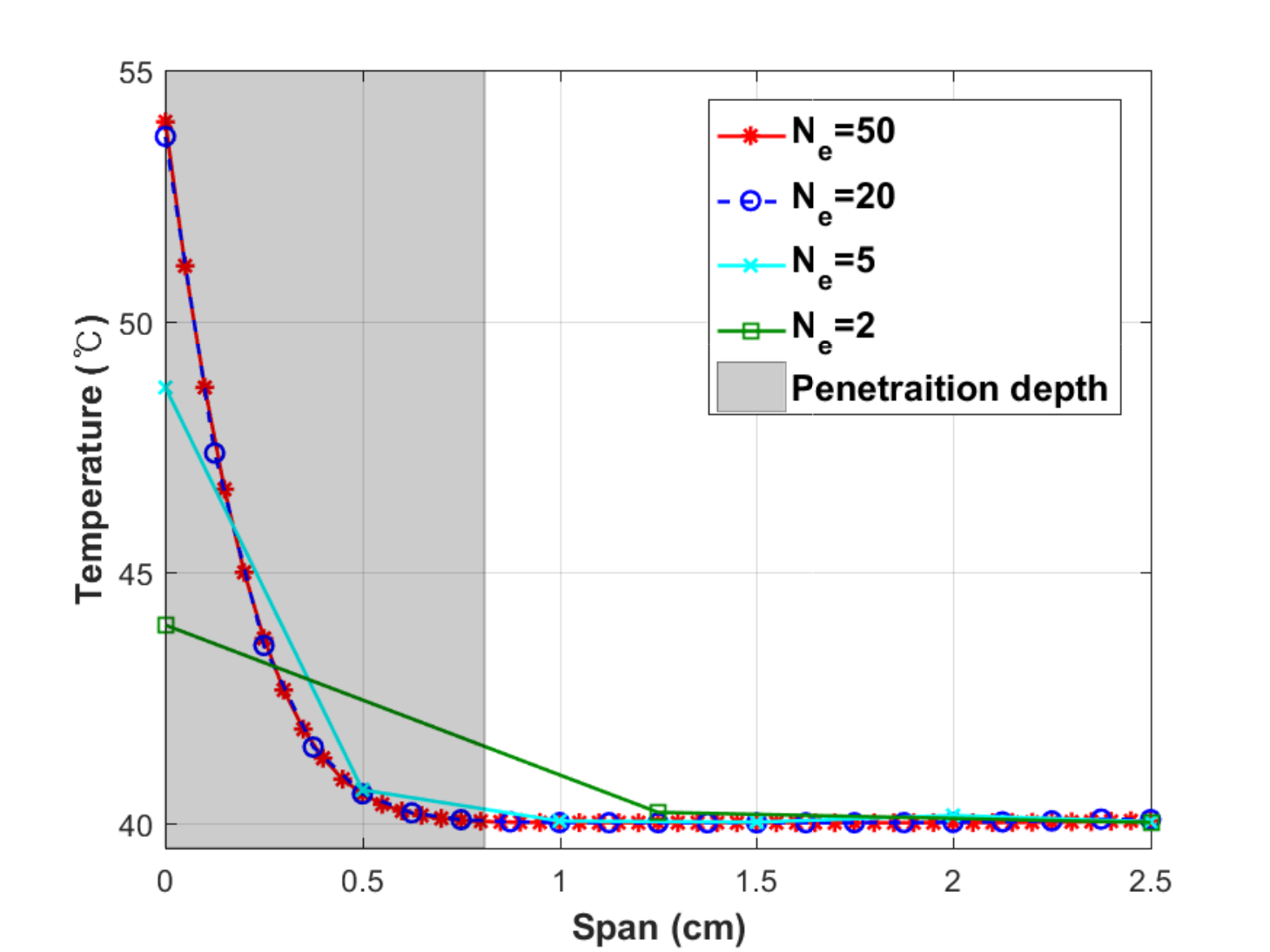}
\caption{  Comparison of FEM results with $\kappa=0.041$, $\Delta t=0.1$, and $\beta=1.0$ at $t=1.0$ to determine the element size in the 1D bar model. Using $N_e=2$ and $N_e=5$ violates the penetration depth condition.  }
    \label{TPDepth}
    \end{figure}

To evaluate the proposed algorithm, one-dimensional heat conduction was implemented in a rectangular bar (Figure~\ref{1Dmodel}). 
The total length $L$ of the bar is 2.5 $cm$, and the cross-section perimeter $P$ and area $A$ are 0.8 $cm$ and 0.04 $cm^2$, respectively.
The true, but unknown, heat flux was imposed on a wall at $x=0$.
The bar is attached to an insulated wall at $x=L$.
The other sides of the bar are exposed to convection under steady-state conditions; thereby, we suppose that the convective heat transfer coefficient is constant and uniform over the bar.
The initial temperature of bar is 40 $^\circ C$.
$N_e$ and $N_g$ are respectively 20 and 21 so that the length of each element is 0.125 $cm$. 
The size of the element was determined by the thermal penetration depth condition~\cite{Bathe1996}, which is defined as
     \begin{equation} \label{TPD}
     \delta \approx \gamma\sqrt{\kappa t},
     \end{equation}
where $\kappa$ is the thermal diffusivity $\left(cm^2/sec\right)$, and the value of $\gamma$ depends on the type of boundary condition and, in this example, is 4.
$\delta$ is the penetration depth that the distance to the surface at which the temperature decreases to 0.5$\%$ of the surface temperature.
To obtain an accurate solution, a sufficiently fine mesh is required inside the penetration depth, and the results in Figure~\ref{TPDepth} show that it is sufficient when 20 elements are used.
In this convergence study, the Backward Euler method ($\beta = 1.0$) was used which is unconditionally stable, and $\Delta t$ was set to 0.1. 
It is assumed that $\Delta t$ should be chosen between 0.1 and 0.01 due to the experimental environment and equipment in this example.
As described in the Introduction, we will control $\beta$, not $\Delta t$, for stability using the condition given in Section~\ref{Beta}.
The accuracy of inverse problems is based on the measurement errors, but when using a discretization method such as FEM, the modeling errors also must be quantified.
The best approach is to compare with the analytical solution like Alifanov~\cite{Alifanov1989}, but we indirectly validated by the convergence study to cover a two-dimensional problem as well.

This section uses one sensor, which is attached to node 2. 
The Gaussian noise $\sigma$ was chosen as a measurement error in the single interior temperature history obtained by the sensor.
The total sensing time $t_M$ to obtain the measurement data is 27 $sec$.
The time, $M_s$, to calculate $\xi(\alpha)$ was set to 2 $sec$ in all examples.

\subsubsection{Verification: $\alpha$ selection}   \label{Verify_a}
\begin{figure}
\centering 
    \begin{subfigure}{0.8\linewidth}
\includegraphics[width=\linewidth]{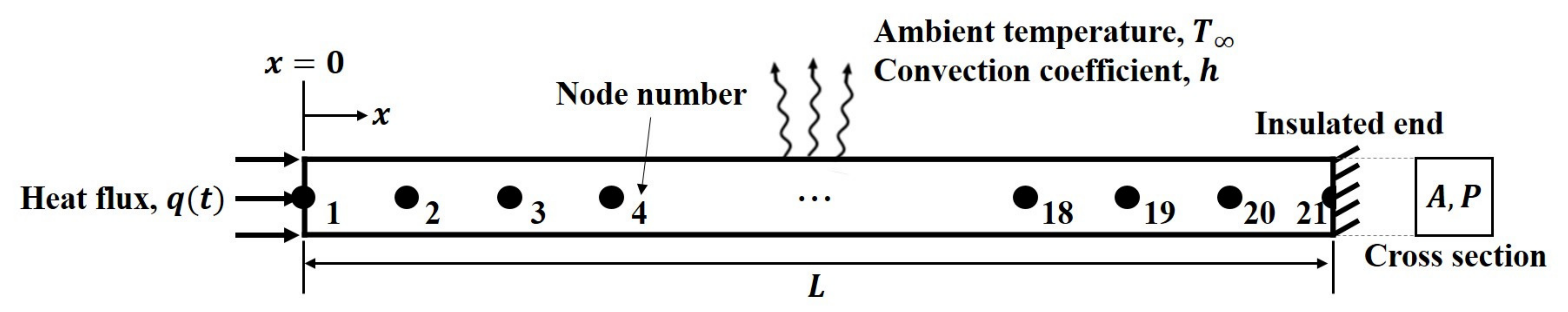}
    \caption{}
    \label{1Dmodel}
    \end{subfigure} \par\bigskip       
    \begin{subfigure}{0.3\linewidth}
\includegraphics[width=\linewidth]{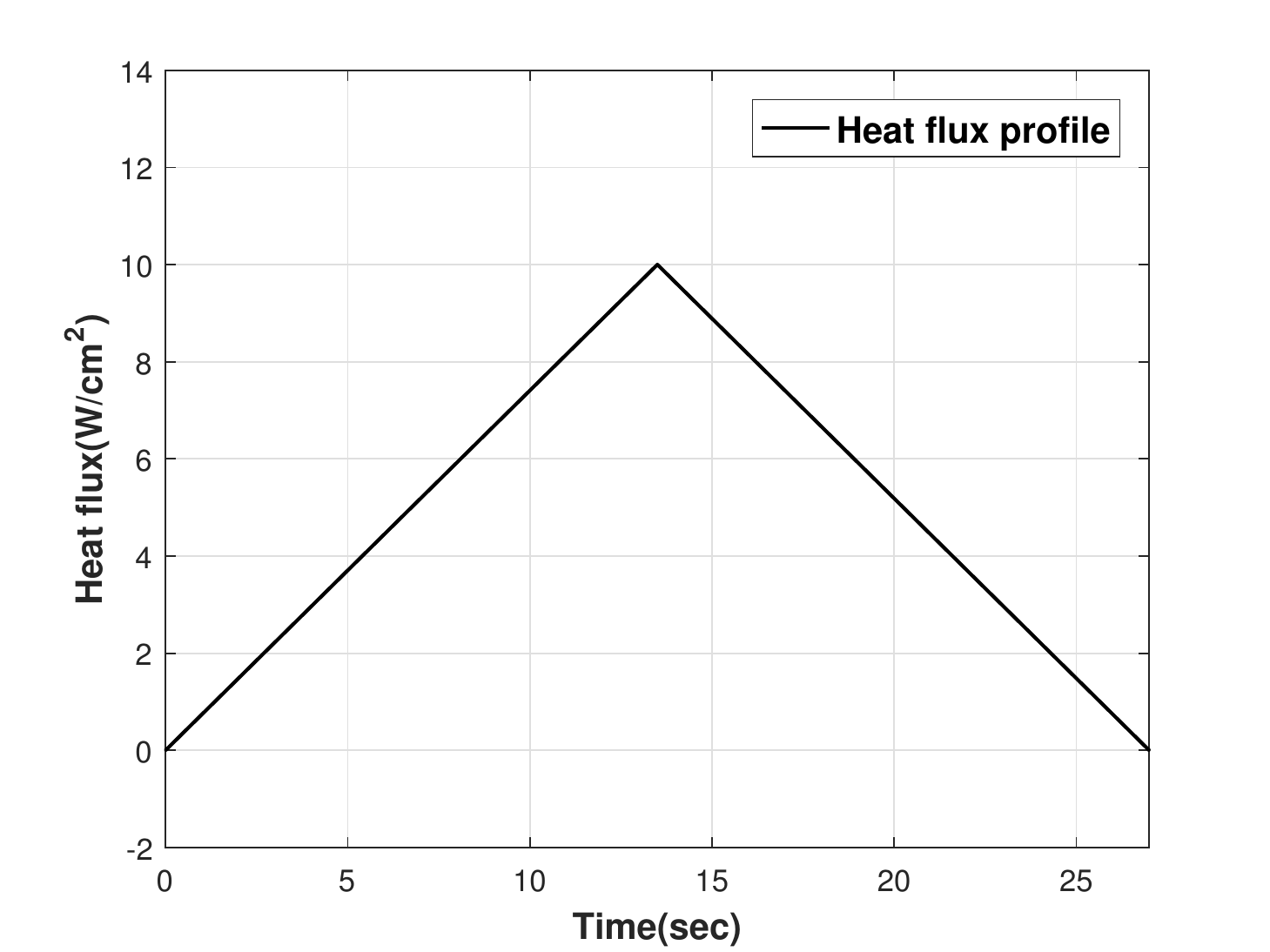}
    \caption{}
    \label{ProfileT}
    \end{subfigure} \hspace{0.5cm}
    \begin{subfigure}{0.3\linewidth}
\includegraphics[width=\linewidth]{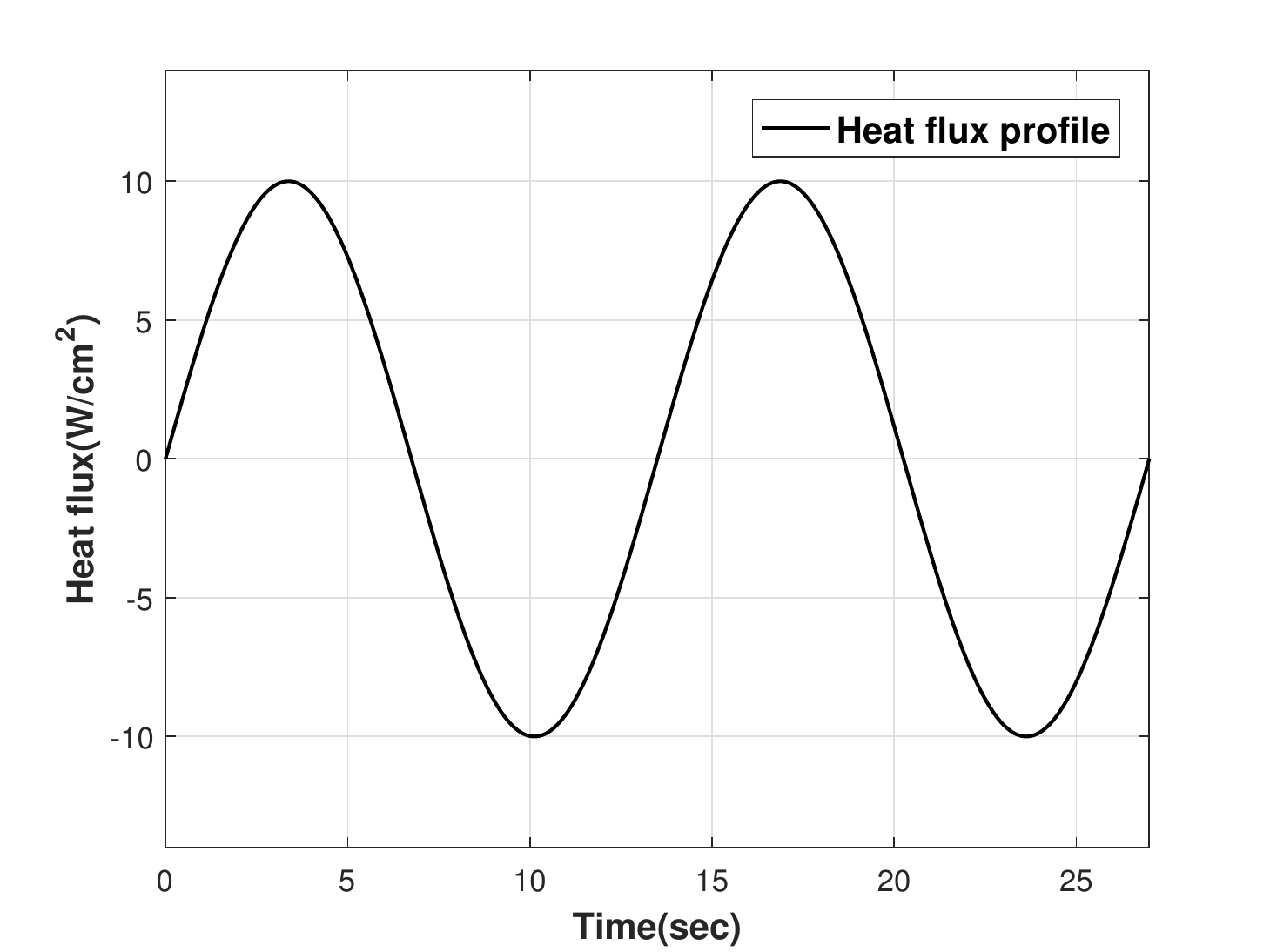}
    \caption{}
    \label{ProfileS}
    \end{subfigure}
\caption{(a) 1D bar finite element model description, and (b) triangular or (c) sinusoidal heat flux profiles imposed on node 1.}
    \label{Profile}
    \end{figure}

\begin{figure}
\centering        
    \begin{subfigure}{0.32\linewidth}
\includegraphics[width=\linewidth]{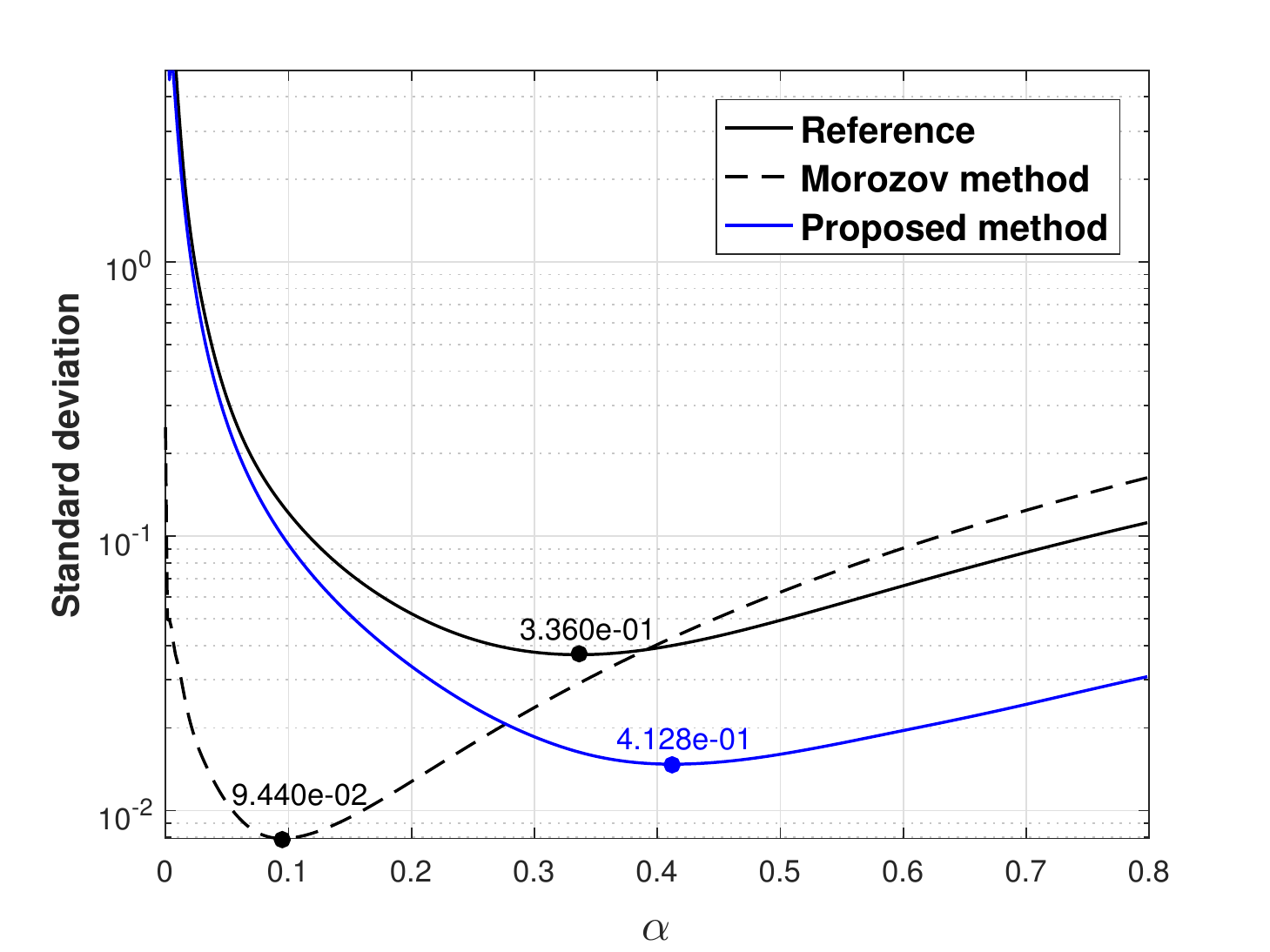}
    \caption{}
    \label{1Dtot006}
    \end{subfigure}
    \begin{subfigure}{0.32\linewidth}
\includegraphics[width=\linewidth]{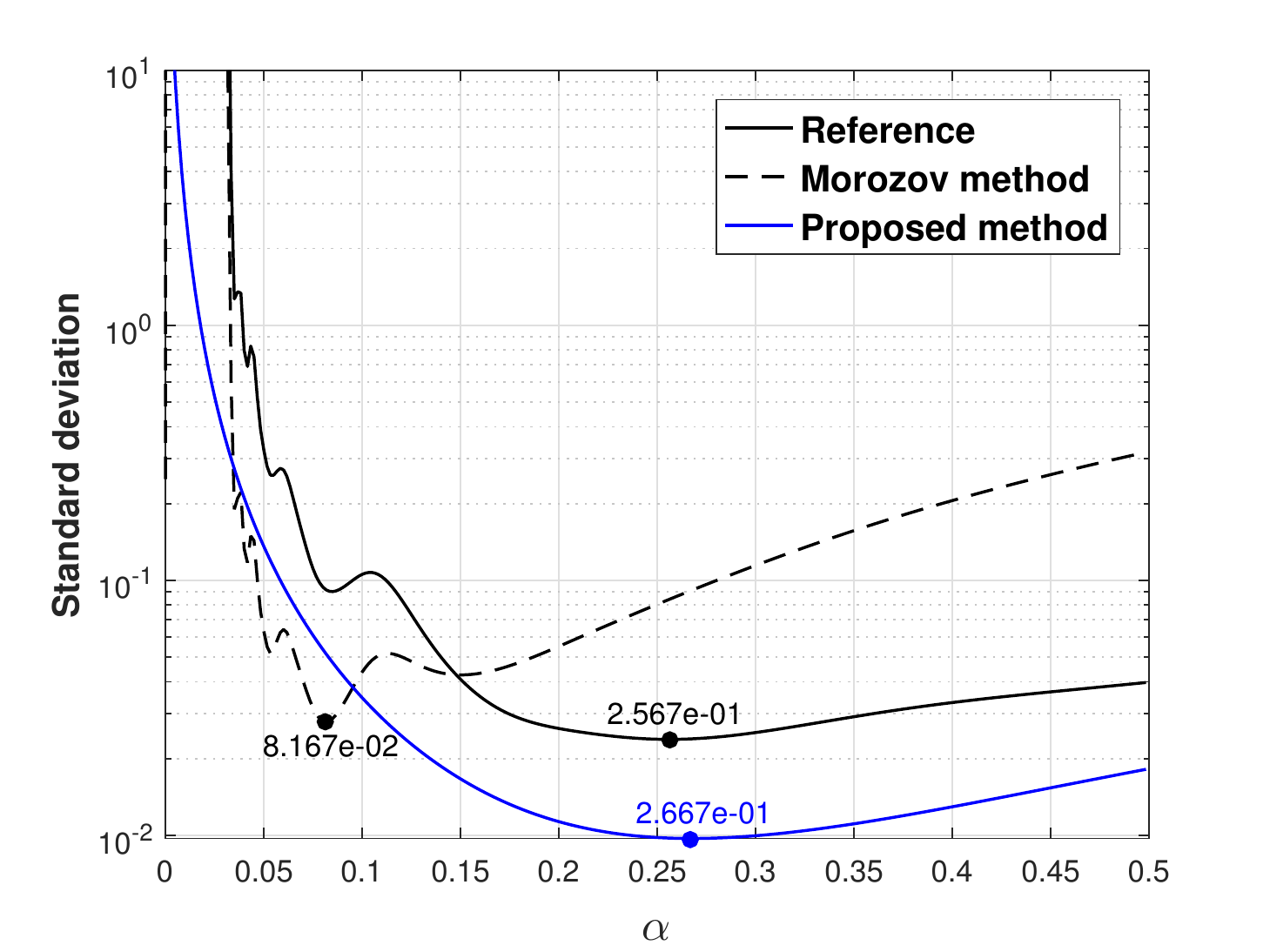}
    \caption{}
    \label{1Dtot006}
    \end{subfigure}
    \begin{subfigure}{0.32\linewidth}
\includegraphics[width=\linewidth]{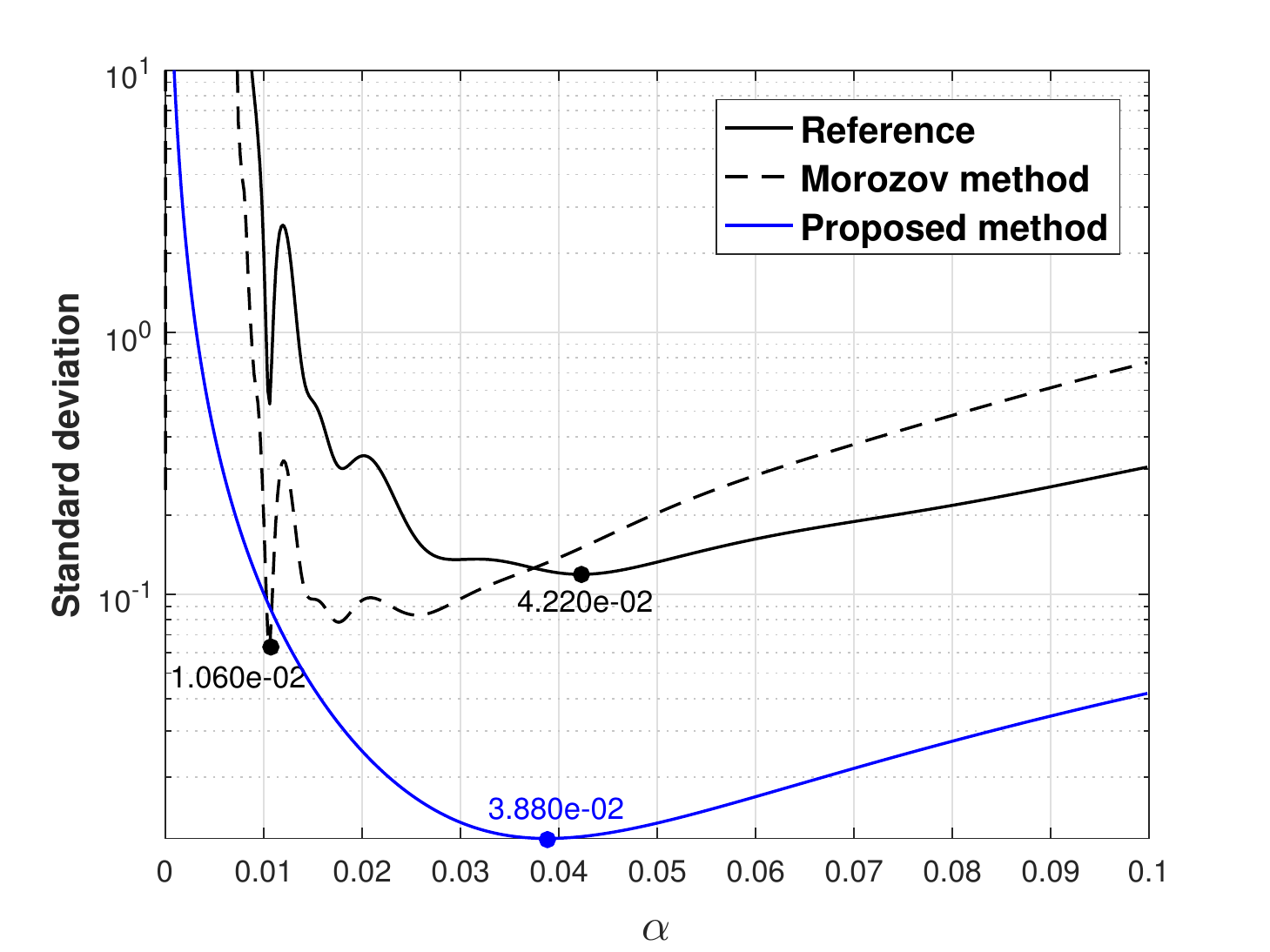}
    \caption{}
    \label{1Dtot006_t}
    \end{subfigure}
\caption{Standard deviation error results of the 1D bar problem when $\alpha$ is varied. 
        Each figure is a result of estimating a triangular heat flux profile for three types of materials: (a) Silicon, (b) Carbon-carbon composite, (c) Stainless steel.
         The average elapsed times to estimate $\alpha$ by the Morozov and the proposed methods are 27.28 and 2.84 $sec$, respectively.}
    \label{Ex1_tot_tri}
    \end{figure}

First, we evaluate the proposed method for selecting the optimal $\alpha$ when $\beta = 1$.
Tests were conducted on two types of heat flux profiles and three types of materials, which are shown in the Figures~\ref{ProfileT},~\ref{ProfileS}, and Table~\ref{Mproperty}.
The optimal $\alpha$ obtained by utilizing the standard deviation in Eq.~\ref{Terror} is used as the reference, and the Morozov discrepancy principle in Eq.~\ref{Morozov} is considered in this verification study.
Equation~\ref{Terror} requires the true heat flux profile, which is unknown in real experimental environments.
Only in numerical studies can the reference of $\alpha$ be calculated using Eq.~\ref{Terror}.
$\Delta q_{max}$ for the triangular and the sinusoidal profiles are 0.0370 and 0.2327 $W/cm^2$, respectively.
$\sigma$ is 0.5, and the sampling time interval $\Delta t$ is 0.01 $sec$.
Figure~\ref{Ex1_tot_tri} shows the standard deviation results from variation in $\alpha$, and also the optimal $\alpha$ obtained by the Morozov and the proposed ridge estimator when $\beta=1.0$.
Herein, the results show that the proposed estimator provides relatively closer results to the reference $\alpha$ than the Morozov method.
In particular, the proposed estimator requires approximately 2 $sec$ to find the optimal $\alpha$, which is 10 times faster than the Morozov method.
The case of sinusoidal heat flux profile was also studied and its results were comparable to the triangular case in the Figure~\ref{Ex1_tot_tri}.

\subsubsection{Verification: $\beta$ selection}  \label{Verify_b}
\begin{figure}
\centering        
    \begin{subfigure}{0.4\linewidth}
\includegraphics[width=\linewidth]{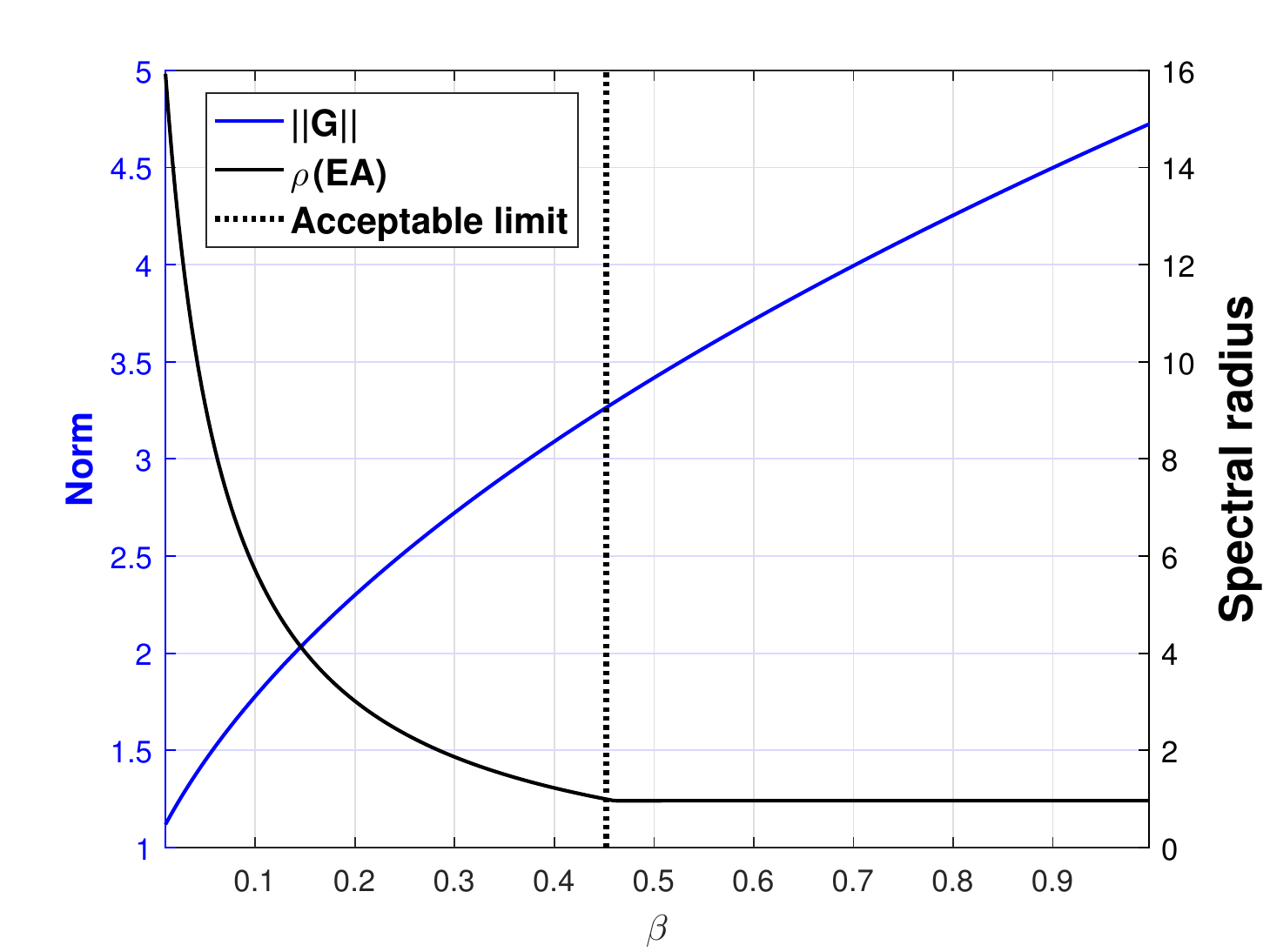}
    \caption{}
    \label{1Dbeta0}
    \end{subfigure} 
    \begin{subfigure}{0.4\linewidth}
\includegraphics[width=\linewidth]{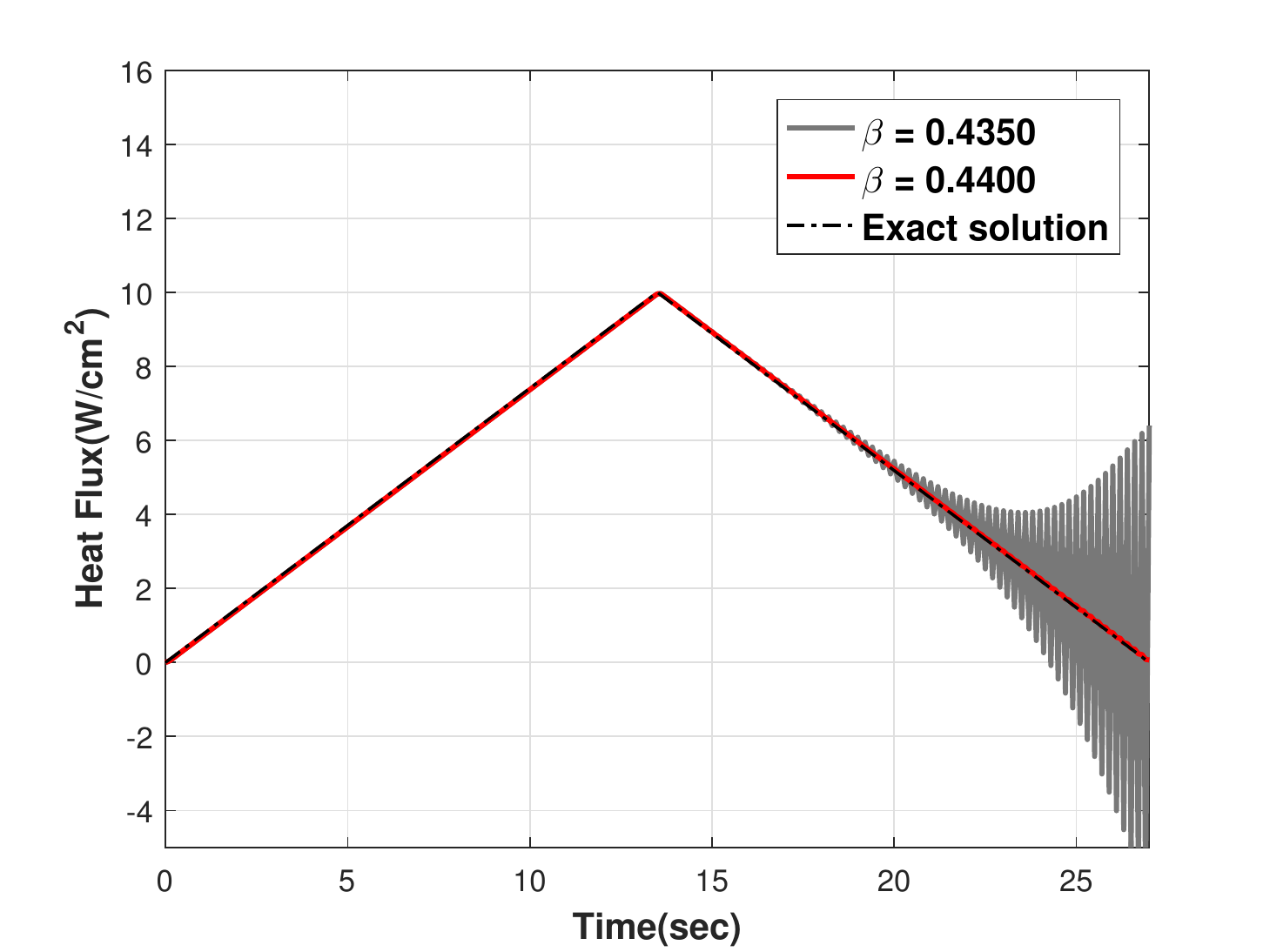}
    \caption{}
    \label{1Dq0}
    \end{subfigure}
\caption{(a) Spectral radius and the norm of the gain coefficient in the 1D bar problem for various values of $\beta$ ($0\le\beta\le1$). 
         The acceptable limit represents the under-bound of the $\beta$ selection to ensure stability.
         (b) Resulting plots of the inverse solutions for validating the stability limit.
         $\sigma$ is set to zero to observe only the initial error propagation.
         The sensor position is node 1 ($\kappa=0.797$, $\Delta t=0.100$, $\alpha=0.0$).}
    \label{Ex1_beta0}
    \end{figure}

In this section, we verify the stability condition studied in Section~\ref{Beta}.
Figure~\ref{Ex1_beta0} shows the norm of the gain coefficient matrix and the spectral radius of $\mb{E}_\alpha\mb{A}$ from variation in $\beta$.
The acceptable limit was drawn where the spectral radius of $\mb{E}_\alpha\mb{A}$ becomes less than 1.
Therefore, $\beta$ with the smallest gain coefficient is selected based on the acceptable limit.
The value of $\beta$ chosen by this condition is 0.44, and the solid red line in Figure~\ref{1Dq0} is the result when $\beta=0.44$.
If $\beta$ is not selected in this manner, the inverse solution may diverge, as shown in Figure~\ref{1Dq0}.

\subsubsection{Performance of the algorithm}  \label{perform1D}
In the previous sections, we evaluated the performance of the proposed method by separately selecting $\alpha$ and $\beta$.
In this problem, we validate the performance of the hybrid parameter selection algorithm by setting two different conditions.
The reference of the inverse solution is calculated by the reference of $\alpha$ and assuming $\beta = 1$.
The Morozov discrepancy principle is also obtained by assuming $\beta = 1$.
However, as shown in Algorithm~\ref{Algorithm}, the hybrid parameters $\alpha$ and $\beta$ are selected by iterative processes.
The errors of the inverse solutions (temperature and heat flux) are calculated by the following formula:
     \begin{subequations}
     \begin{gather}
     \textbf{Error of temp.} = \frac{1}{M}\sum_{m=1}^{M}\left[ \tilde{\mb{T}}^m - \mb{T}_{\text{method}}^m \right]^T\left[ \tilde{\mb{T}}^m - \mb{T}_{\text{method}}^m \right],\\
     \textbf{Error of heat flux} = \frac{1}{M}\sum_{m=1}^{M}\left[ \tilde{\mb{q}}^m - \mb{q}_{\text{method}}^m \right]^T\left[ \tilde{\mb{q}}^m - \mb{q}_{\text{method}}^m \right],
     \end{gather}
     \end{subequations}
where $\tilde{\mb{T}}^m$ and $\tilde{\mb{q}}^m$ are the exact temperature and heat flux vectors, respectively. 
$\mb{T}_{\text{method}}^m$ and $\mb{q}_{\text{method}}^m$ are obtained by each method individually: the reference, Morozov method, and hybrid parameter selection method.

\begin{figure}
\centering        
    \begin{subfigure}{0.4\linewidth}
\includegraphics[width=\linewidth]{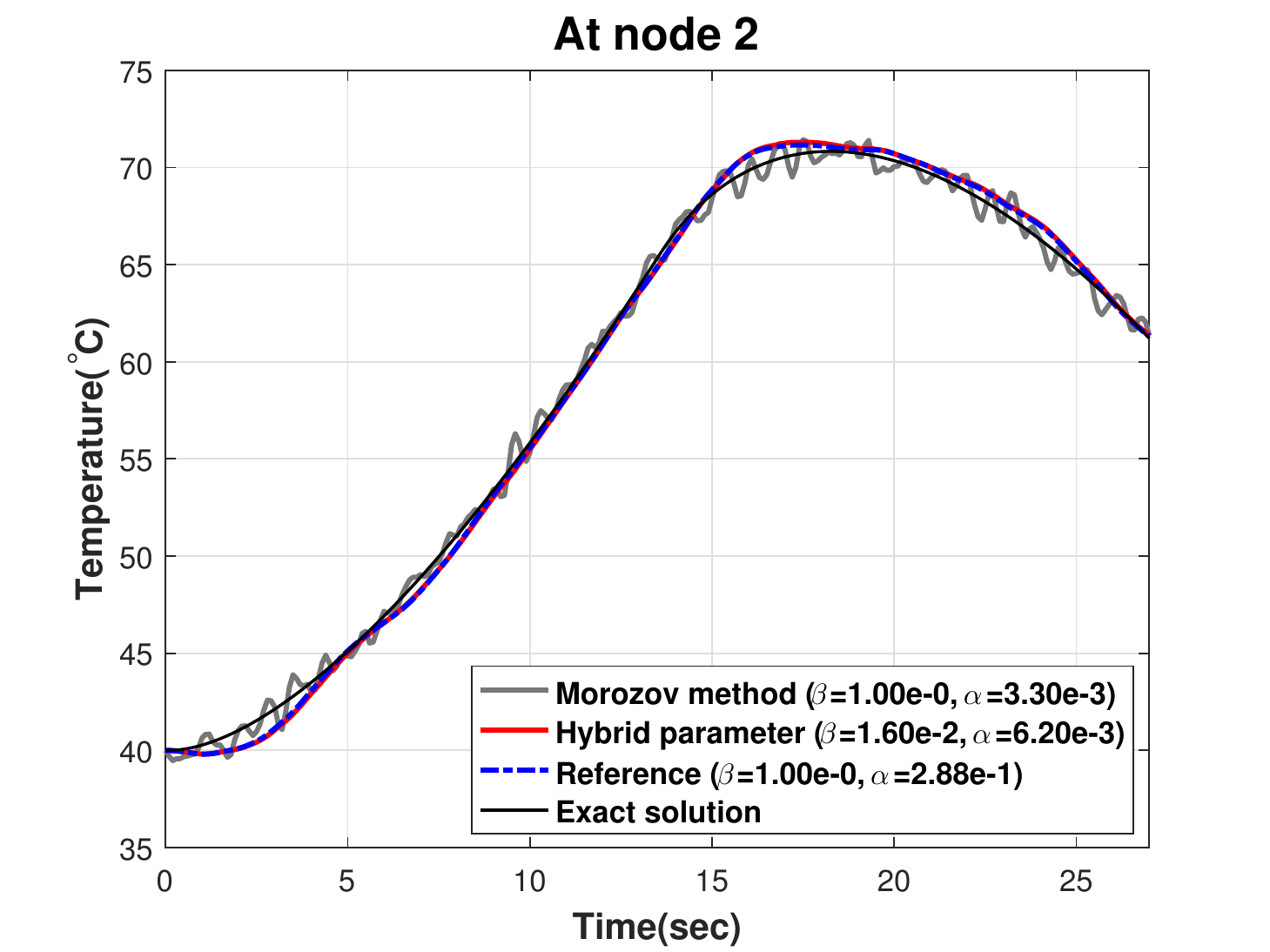}
    \caption{}
    \label{1Dhybrid_tri1}
    \end{subfigure}
    \begin{subfigure}{0.4\linewidth}
\includegraphics[width=\linewidth]{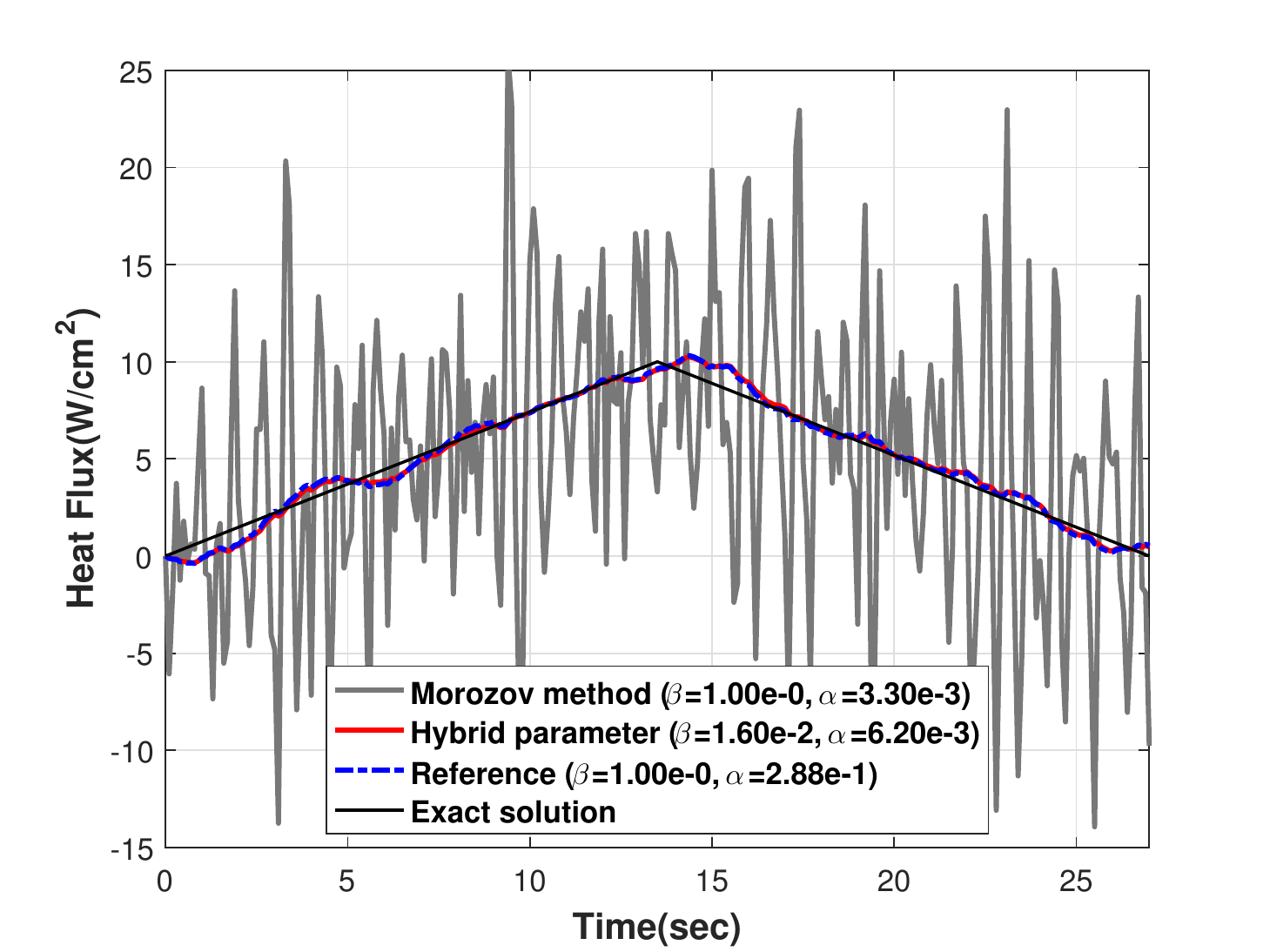}
    \caption{}
    \label{1Dhybrid_sine1}
    \end{subfigure}
\caption{Estimated temperature (left) and heat flux (right) plots using the Morozov and the hybrid parameter methods for the triangular input cases in the 1D bar problem.
         Reference is the result using the optimal regularization parameter $\alpha$ obtained by the exact solution.
         The sensor position and thermal diffusivity are presented in Table~\ref{CompResults2}.}
    \label{1Dhybrid1}
    \end{figure}

Figure~\ref{1Dhybrid1} shows the estimated temperature and heat flux plots.
The standard deviation errors and computational cost are shown in Table~\ref{CompResults2}.
The case of sinusoidal heat flux profile was also studied and its results are shown in Table~\ref{CompResults2}.
The hybrid parameter selection method reduced the optimal value of $\alpha$ significantly, and
the standard deviation of the heat flux calculated by the proposed algorithm is similar to that of the reference.
Table~\ref{CompResults2} indicates that the Morozov method yields a relatively large error for heat flux, which is because this method is a temperature-based estimator.


\subsection{Two-dimensional problem} \label{EX2D}
A two-dimensional heat conduction problem was implemented in a rectangular plate to demonstrate the accuracy and efficiency of the hybrid parameter selection algorithm.
The thickness of the plate is thin enough to ignore the heat conduction in the through-plane, and
the description of the model is shown in Figure~\ref{Ex2}. 
This study considers natural convection assuming that the plate is exposed to ambient room air without an external source of motion, and the initial equilibrium temperature is 40 $^\circ C$.
The size of the plate is (5$cm$ $\times$ 5$cm$ $\times$ 0.01$cm$), featuring 400 quad elements (i.e., 2D elements) and 441 nodes.
The size of the element was determined in the same way as in Section 4.1.
It is assumed that the measurement time interval must be 0.05 (i.e., $\Delta t = 0.05$).
The material density $\rho$, thermal conductivity $k$, specific heat capacity $c_v$, and convective heat transfer coefficient $h$
are 1.33 $g/cm^3$, 2.33 $W/cm ^\circ C$, 0.940 $J/g ^\circ C$, and 0.0005 $W/cm^2$$^\circ C$, respectively.
The time, $M_s$, to calculate $\xi(\alpha)$ was set to 2 $sec$ in all examples, and $\xi(\alpha)$ was used as a fitted curve.
The reference and the error of the inverse solution were calculated in the same manner as in Section~\ref{EX1D}.
In the following sections, two cases were considered: single-unknown and multi-unknown heat fluxes.

\subsubsection{Two-dimensional problem with single-unknown heat flux}  \label{perform2D}
\begin{figure}
\centering        
    \begin{subfigure}{0.37\linewidth}
\includegraphics[width=\linewidth]{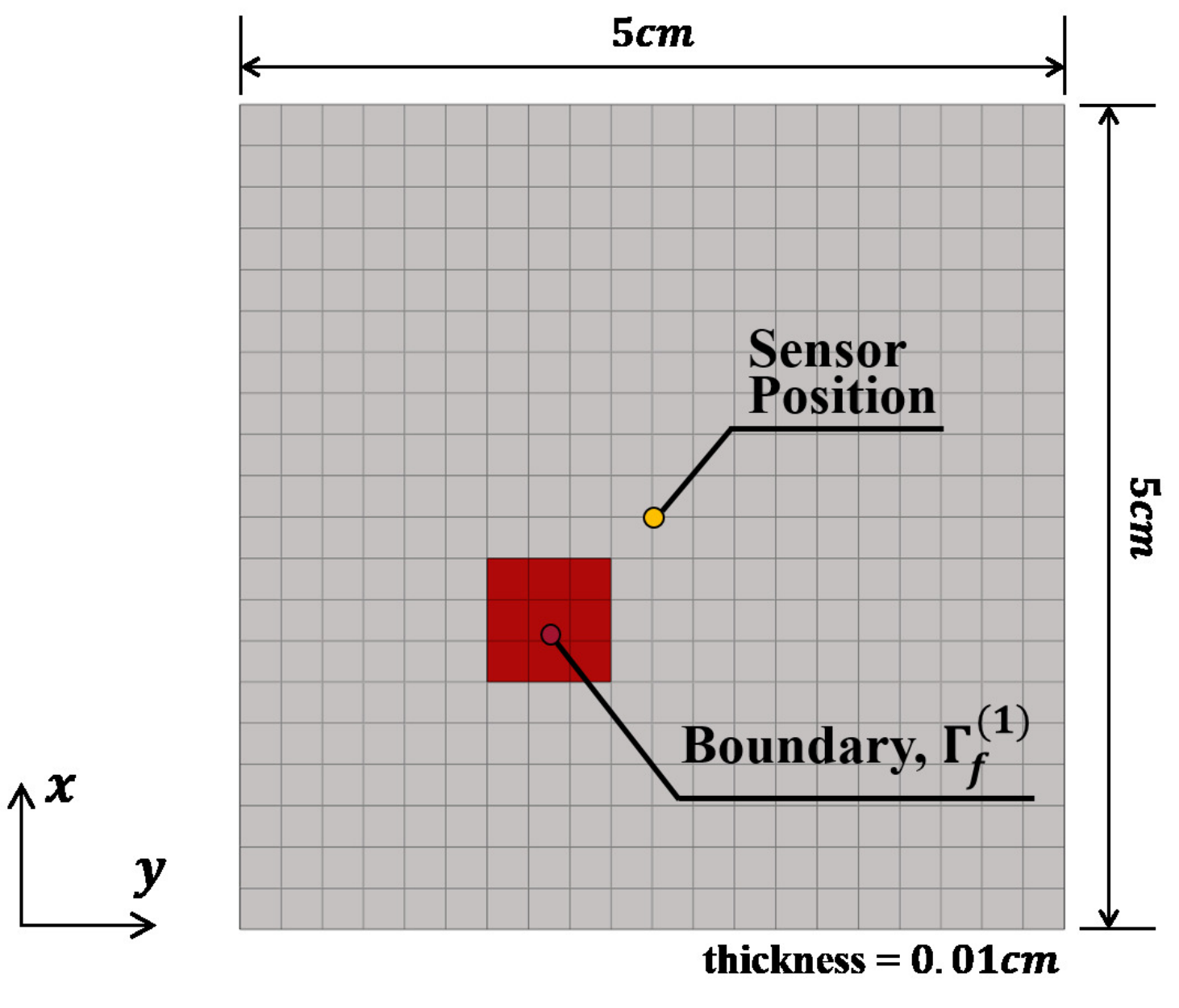}
    \caption{}
    \label{Ex2_model}
    \end{subfigure}  \hspace{1.5cm}
    \begin{subfigure}{0.35\linewidth}
\includegraphics[width=\linewidth]{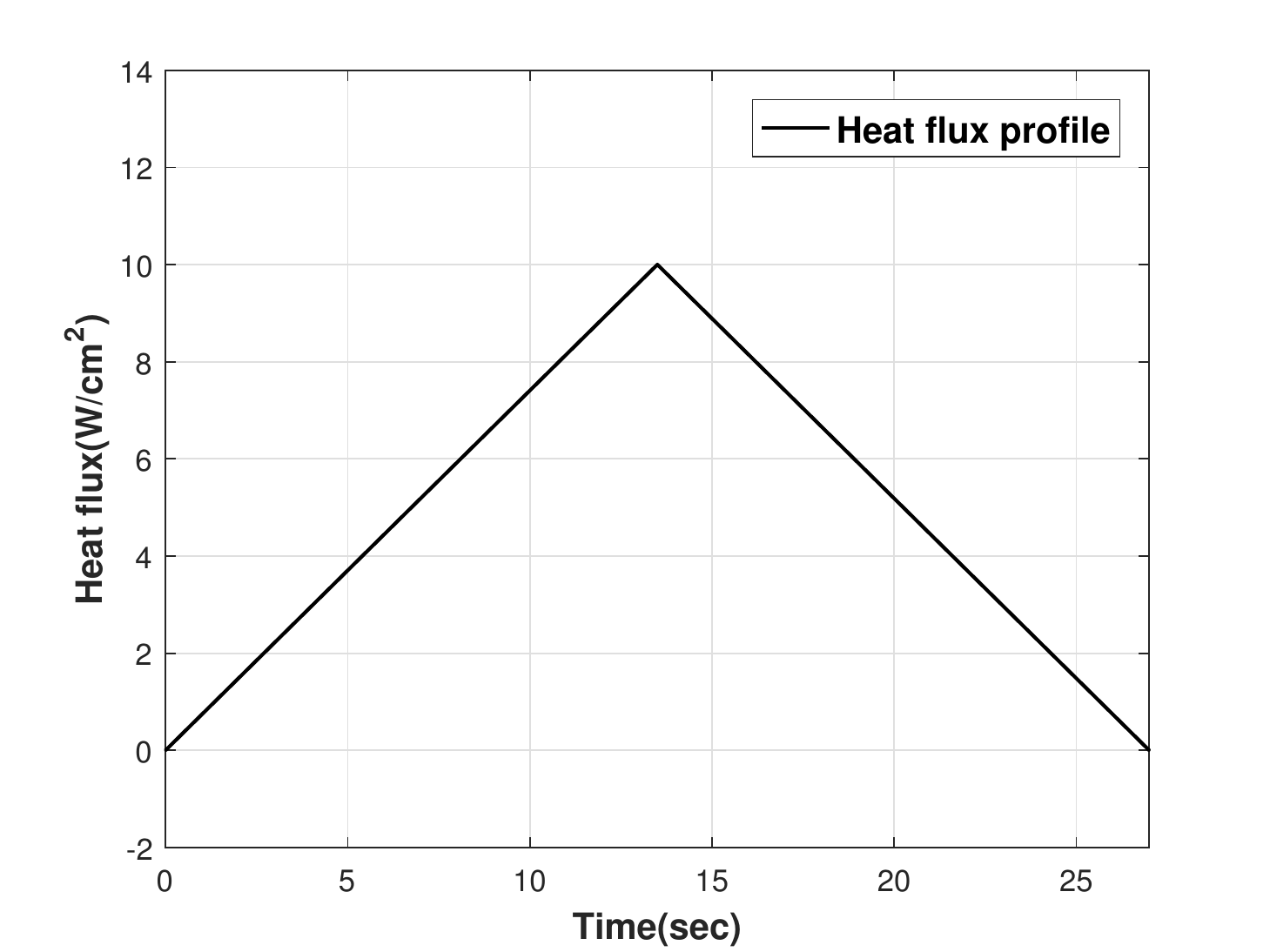}
    \caption{}
    \label{2DProfile}
    \end{subfigure}
\caption{(a) 2D finite element model and sensor position, and (b) heat flux is imposed on the boundary $\Gamma_f^{(1)}$.}
    \label{Ex2}
    \end{figure}

\begin{figure}
\centering    
    \begin{subfigure}{0.57\linewidth}
\includegraphics[width=\linewidth]{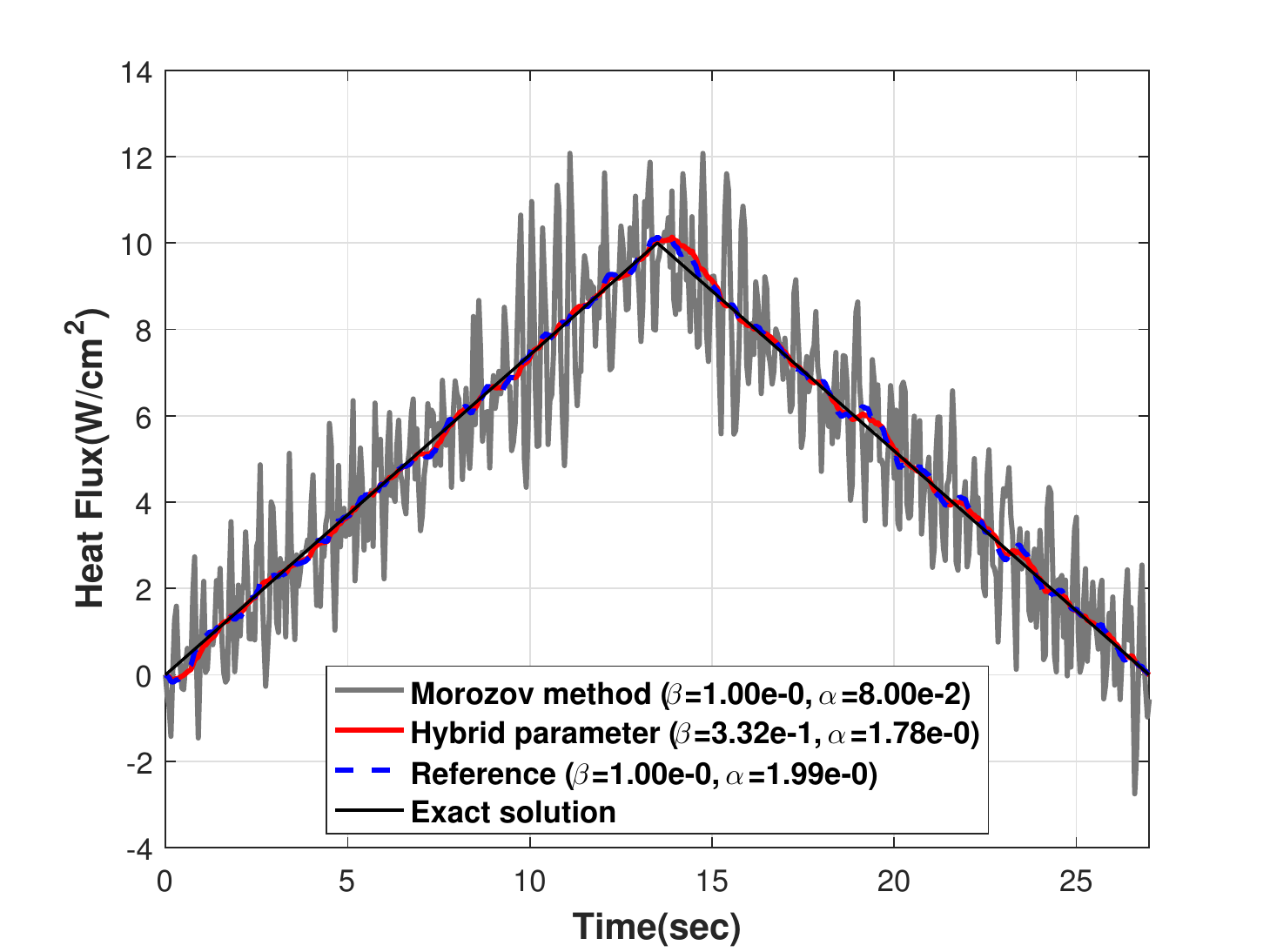}
    \caption{}
    \end{subfigure}
\caption{Estimated heat flux plots for the triangular input case in the 2D plate problem.
         The average elapsed times to execute the Morozov and the hybrid parameter algorithms are 452.26 and 71.79 $sec$, respectively ($\kappa=1.863$, $\sigma=0.5$, $\Delta t=0.05$).}
    \label{Ex2_hybrid1}
    \end{figure}

The first case is to estimate a single-unknown boundary heat flux.
The locations of the imposed heat flux and the attached sensor are shown in Figure~\ref{Ex2_model}. 
The profile of the imposed heat flux is shown in Figure~\ref{2DProfile}.
The result is shown in Figure~\ref{Ex2_hybrid1}, and the details of the results are listed in Table~\ref{CompEx2_hybrid1}.
The numerical results clearly demonstrate that the proposed hybrid parameter selection algorithm yields better accuracy than the Morozov method when estimating the heat flux.
In addition, it requires only 71.79 $sec$ to compute the two parameters, which is much faster than the Morozov method (452.26 $sec$).

\subsubsection{Two-dimensional problem with multi-unknown heat fluxes}  \label{perform2Dmulti}
\begin{figure}
    \centering
    \begin{subfigure}[t]{0.3\linewidth}
\includegraphics[width=\linewidth]{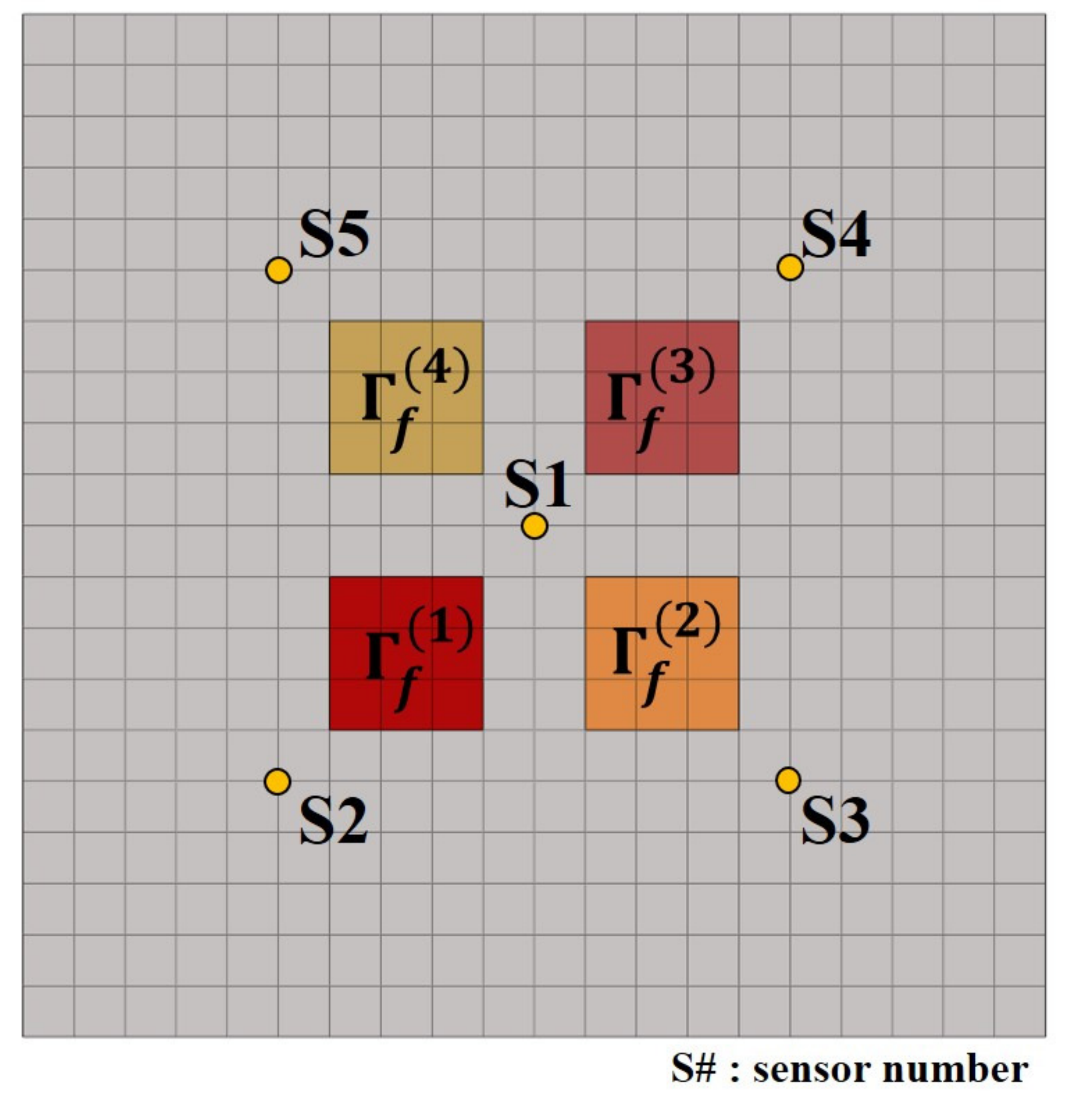}
    \caption{}
    \label{Ex2_model2}
    \end{subfigure}  \hspace{1cm}
    \begin{subfigure}[t]{0.5\linewidth}
\includegraphics[width=\linewidth]{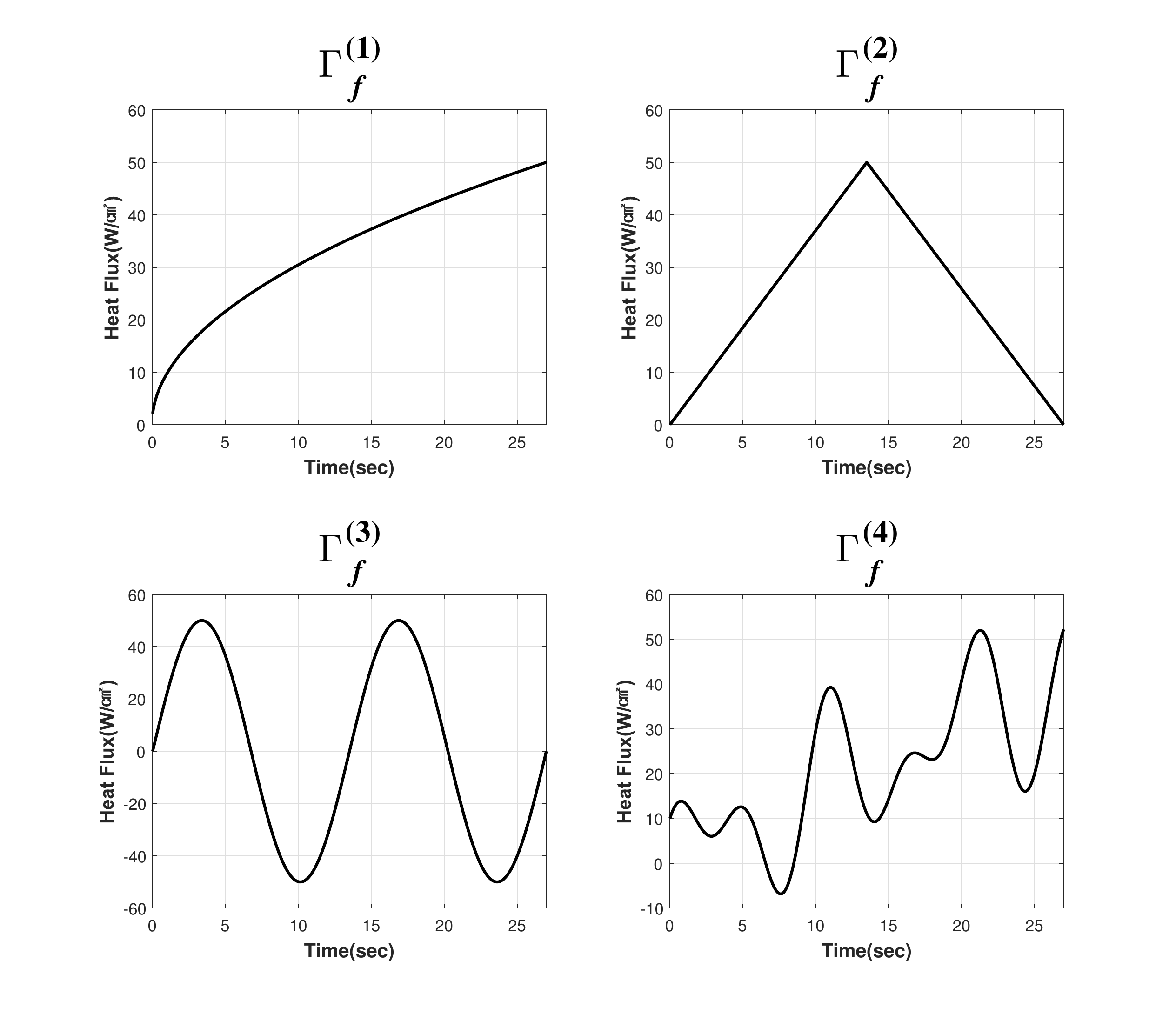}
    \caption{}
    \label{2DProfileMulti}
    \end{subfigure}
\caption{(a) Description of the attached sensors and the imposed boundary conditions. 
         (b) The different types of the heat flux profiles are imposed on the boundary ($\Gamma_f^{(1)}$: square root, $\Gamma_f^{(2)}$: triangular, $\Gamma_f^{(3)}$: sinusoidal, $\Gamma_f^{(4)}$: arbitrary profiles). }
    \label{Ex2_2}
    \end{figure}

\begin{figure}
\centering    
    \begin{subfigure}{0.4\linewidth}
\includegraphics[width=\linewidth]{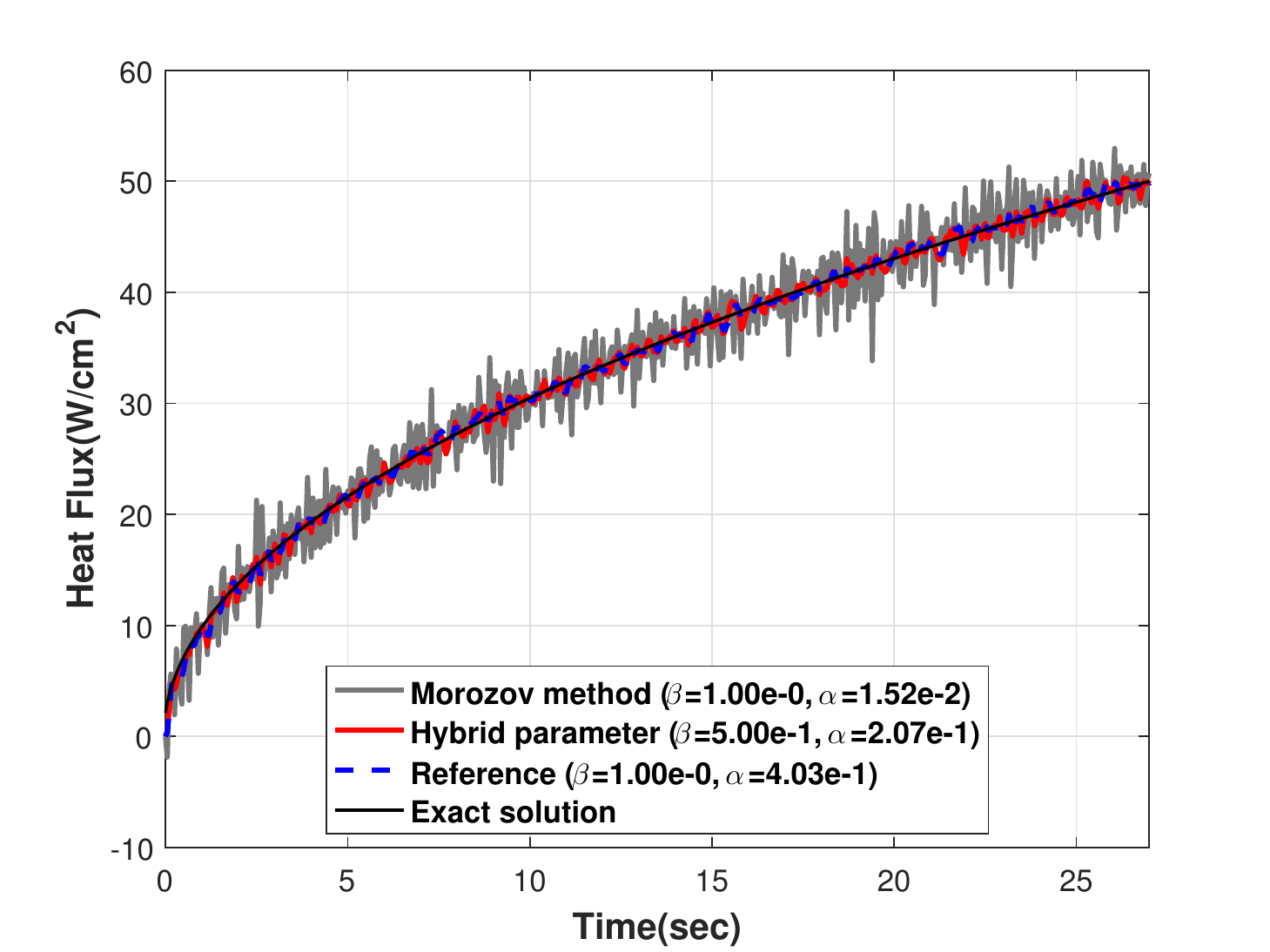}
    \caption{$\Gamma_f^{(1)}$}
    \label{Ex2_multi1}
    \end{subfigure} 
    \begin{subfigure}{0.4\linewidth}
\includegraphics[width=\linewidth]{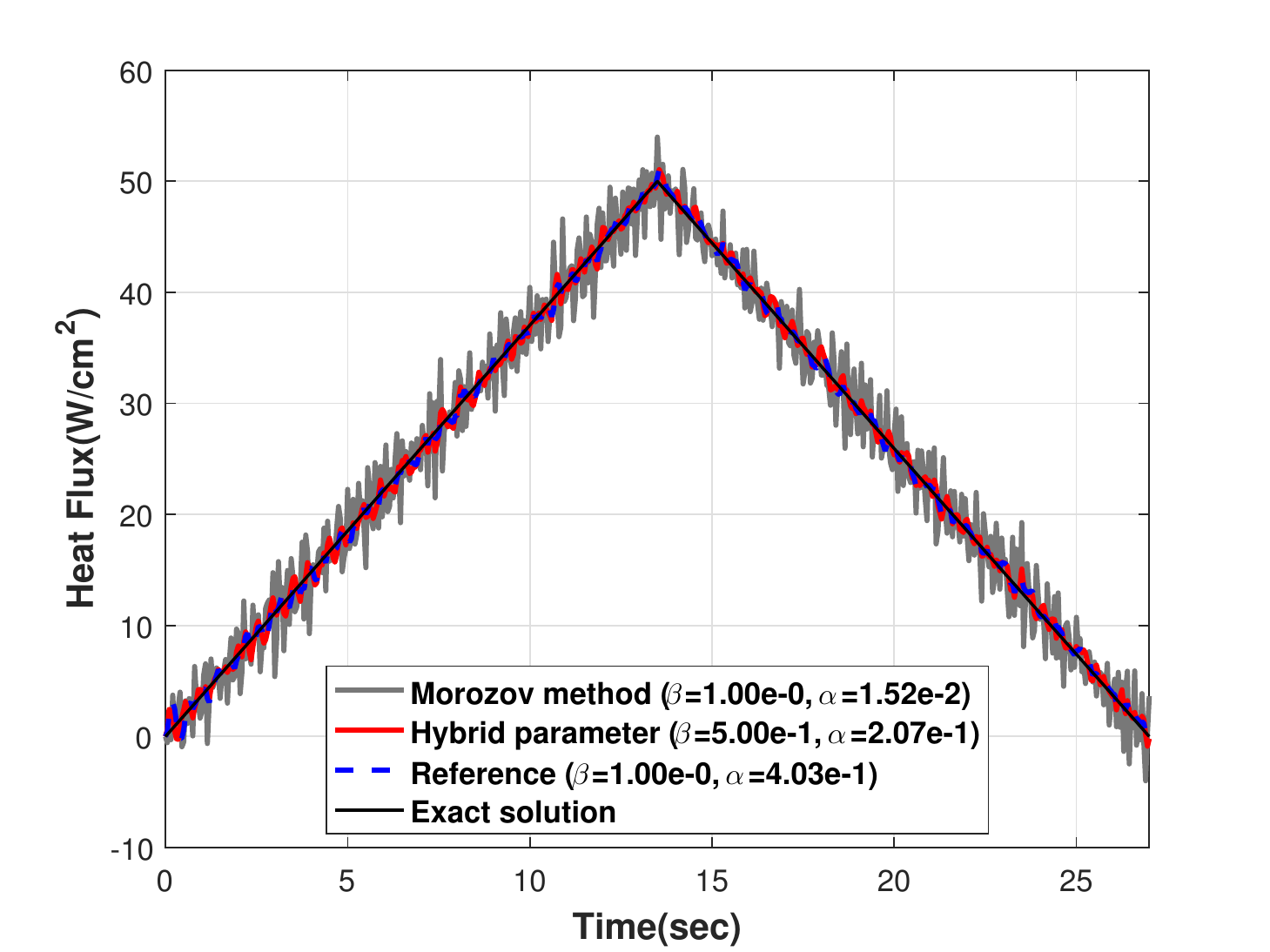}
    \caption{$\Gamma_f^{(2)}$}
    \label{Ex2_multi2}
    \end{subfigure}
    \begin{subfigure}{0.4\linewidth}
\includegraphics[width=\linewidth]{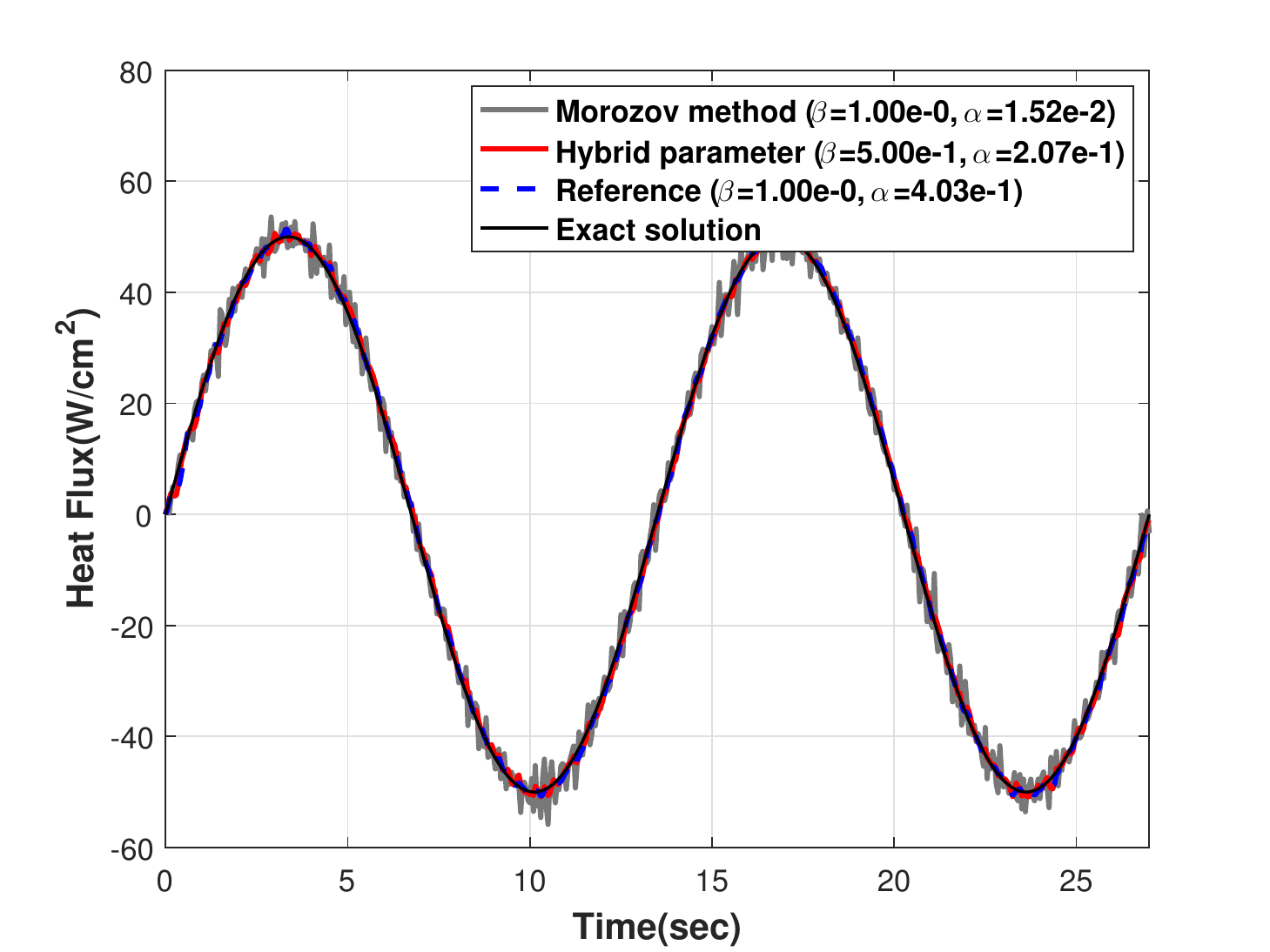}
    \caption{$\Gamma_f^{(3)}$}
    \label{Ex2_multi3}
    \end{subfigure}
    \begin{subfigure}{0.4\linewidth}
\includegraphics[width=\linewidth]{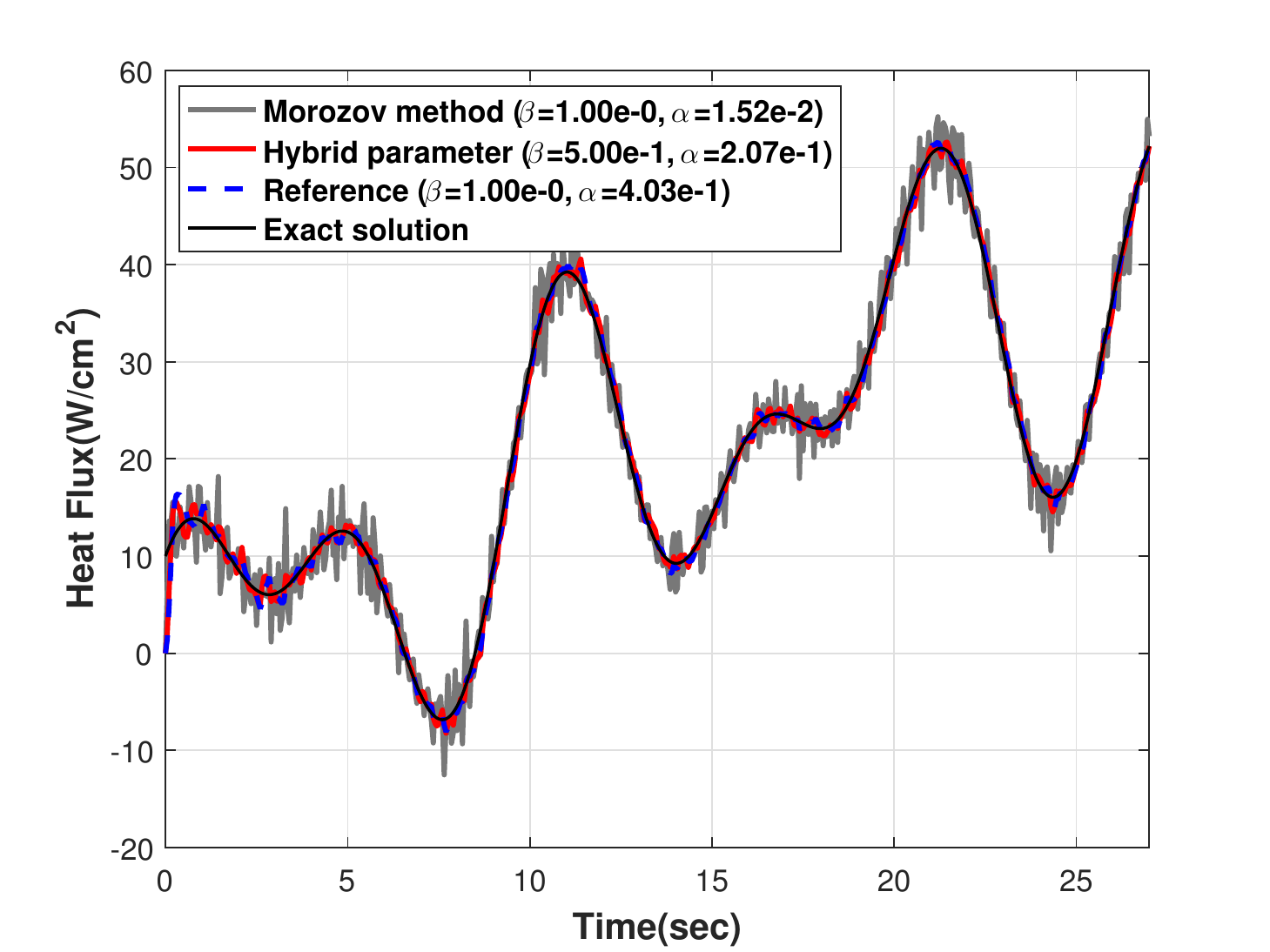}
    \caption{$\Gamma_f^{(4)}$}
    \label{Ex2_multi4}
    \end{subfigure}
\caption{Results for the multiple-unknown heat flux estimation in the 2D plate problem.
         The average elapsed times to execute the Morozov and the hybrid parameter algorithms are 475.43 and 73.56 $sec$, respectively ($\kappa=1.863$, $\sigma=0.5$, $\Delta t=0.05$).}
    \label{Ex2_multi}
    \end{figure}

The second case is to estimate multi-unknown boundary heat fluxes.
The description of the model and boundary conditions is shown in Figure~\ref{Ex2_model2}. 
The profiles of the heat flux imposed for each boundary are shown in Figure~\ref{2DProfileMulti}.
The results are shown in Figure~\ref{Ex2_multi}, and the details of the results are listed in Table~\ref{CompEx2_multi}.
With the corresponding results, the proposed ridge estimator can be applied well to the case with multiple sensors.
In addition, for the problem of estimating multiple heat sources, the hybrid parameter selection method was successfully validated in terms of accuracy and efficiency.

\section{Conclusion}  \label{Conclusion}
In this study, the iterative hybrid parameter ($\alpha$ and $\beta$) selection algorithm was developed.
$\alpha$ is the regularization parameter that gives the bias error to reduce the variance error, and 
$\beta$ is the parameter to reduce the bias error.
Thus, the proposed algorithm can reduce the total error by adjusting the two parameters iteratively, resulting in better stability.
An efficient ridge estimator was first proposed to alleviate the computational burden of the iterative scheme. 
The key idea is to decompose the total error into bias and variance errors, and each error is effectively computed using the maximum heat flux change and prior information, respectively.   
For this purpose, the proposed ridge estimator possesses two features.
One is that it does not require experimental data, which needs to be obtained in advance.
Therefore, the proposed estimator can be used in pre-processing as well, whereas classical methods can only be used in post-processing.
Another feature is that it does not have a local minimum unlike the Morozov method.
Thus, whereas the Morozov method has to specify the range of $\alpha$ (i.e., algorithms designed for-loop), the proposed estimator can be efficiently implemented by using the while-loop in programming.
In the hybrid parameter selection scheme, the amplification matrix of the inverse problem derived by the Tikhonov method was derived here.
$\beta$ was then selected to obtain better accuracy of the inverse solution when the spectral radius of the amplification matrix is less than 1.  
Combining these criteria, a new iterative hybrid parameter selection algorithm was developed to find $\alpha$ and$\beta$, and its performance was demonstrated using numerical examples.
The proposed algorithm can address the multiple heat source estimation problem requiring multi sensor data. 
Although this work focused only on IHCP, the idea could be extended to other problems such as the thermal expansion, piezoelectric, inverse radiation, natural convection, and elastic problems~\cite{Rabczuk2016, Rabczuk2019, Zabaras1997, NMPark2003, KHLee2008, Jaluria2021}.

Note that the deviation error was neglected in the proposed estimator.
However, this error becomes significant when sensors are attached in wrong positions. 
Therefore, sensitivity analysis should be conducted to find adequate sensor positions before using the proposed algorithm.
Furthermore, the proposed estimator cannot be used to estimate a step or constant input case because it is based on the heat flux change due to the computational efficiency.

\section*{Acknowledgements}
This research was funded by National Research Foundation of Korea (NRF-2018R1A1A1A05078730) and MOTIE, Korea Government (grant No. 20001228).


\clearpage

\clearpage
\clearpage

\begin{table}
  \begin{center}    
    \caption{Material properties for one-dimensional numerical examples.}
    \label{Mproperty}
    \renewcommand{\arraystretch}{1.5}
    \begin{tabular}{l c c c} \hline
                                                            &        \MCC{\textbf{Material properties}}   \\\cline{2-4}
                                                            &\MR{\textbf{Silicon}}&\textbf{Carbon-carbon} &\MR{\textbf{Stainless steel}} \\[-7pt]
                                                            &                     &\textbf{composite}     &          \\\hline
      \text{Thermal diffusivity}, $\kappa$ [$cm^2/sec$]     &\text{0.797}    &\text{0.103}           &\text{0.041}   \\
      Material density, $\rho$ [$g/cm^3$]                   &    2.330       &     1.720             &     7.860     \\
      Thermal conductivity, $k$ [$W/cm ^\circ C$]           &    1.300       &     0.125             &     0.162     \\
      Specific heat capacity, $c_v$ [$J/g ^\circ C$]        &    0.700       &     0.700             &     0.500     \\
      Convection coefficient, $h$ [$W/cm^2$$^\circ C$]      &                  \MCC{0.005}                \\
      Ambient temperature, $T_{\infty}$ [$^\circ C$]        &                  \MCC{40.000}               \\
      Initial temperature, $T_0$ [$^\circ C$]               &                  \MCC{40.000}               \\\hline
    \end{tabular}
  \end{center} 
\end{table}

\clearpage

\clearpage

\clearpage
\begin{table}
  \begin{center}    
    \caption{ Detail results for each heat flux profile case in the 1D bar problem. }
    \label{CompResults2}
    \renewcommand{\arraystretch}{1.2}
    \begin{tabular}{c || c c c c c c | c} \toprule[1.2pt]
                             &\MR{\textbf{Method}}& \MR{\alpha}  & \MR{\beta}    &\textbf{Error of}  &\textbf{Error of} &\MR{\textbf{Cost}\text{(sec)}}&\MR{\textbf{Note}}\\
                             &                    &              &               &\textbf{temp.}     &\textbf{heat flux}&                              &                  \\\midrule[1.2pt]
      \textbf{Triangular}    & Reference          &   2.88$e$-1  &    1.00$e$-0  &       4.75$e$-1   &  4.54$e$-1       &  3.13$e$-0   &   sensor attached to node 2,     \\
      \textbf{heatflux}      & Morozov method     &   3.30$e$-3  &    1.00$e$-0  &       5.08$e$-1   &  6.69$e$-0       &  3.38$e$-0   &        $\kappa=0.041$,           \\
      \textbf{profile}       & Hybrid parameter   &   6.20$e$-3  &    1.60$e$-1  &       5.19$e$-1   &  4.48$e$-1       &  9.11$e$-1   &  $\sigma = 0.5$, $\Delta t=0.1.$ \\\hline
      \textbf{Sinusoidal}    & Reference          &   3.47$e$-2  &    1.00$e$-0  &       2.97$e$-1   &  5.81$e$-1       &  2.74$e$+1   &   sensor attached to node 3,     \\
      \textbf{heatflux}      & Morozov method     &   1.13$e$-2  &    1.00$e$-0  &       1.72$e$-1   &  9.21$e$-1       &  2.92$e$+1   &        $\kappa=0.797$,           \\
      \textbf{profile}       & Hybrid parameter   &   1.99$e$-5  &    1.60$e$-1  &       2.86$e$-1   &  5.84$e$-1       &  4.60$e$-0   & $\sigma = 0.5$, $\Delta t=0.01.$ \\\bottomrule[1.2pt]
    \end{tabular}
  \end{center}  
\end{table}

\clearpage
\begin{table}
  \begin{center}    
    \caption{ Detail results for the case of a single-unknown heat flux profile estimation in the 2D plate problem. }
    \label{CompEx2_hybrid1}
    \renewcommand{\arraystretch}{1.2}
    \begin{tabular}{c | c c c c c | c} \toprule[1.2pt]
    \MR{\textbf{Method}}& \MR{\alpha}  & \MR{\beta}    &\textbf{Error of}  &\textbf{Error of} &\MR{\textbf{Cost}\text{(sec)}}&\MR{\textbf{Note}} \\
                        &              &               &\textbf{temp.}     &\textbf{heat flux}&                              &                   \\\midrule[1.2pt]
     Reference          &   1.99$e$-0  &    1.00$e$-0  &       3.77$e$-1   &  1.53$e$-1       &  4.48$e$+2   &     \MR{\kappa=1.863,}            \\
     Morozov method     &   8.00$e$-2  &    1.00$e$-0  &       1.24$e$-1   &  1.22$e$+0       &  4.52$e$+2   &\MR{\sigma = 0.5, \: \Delta t=0.05}\\
     Hybrid parameter   &   1.78$e$-0  &    3.32$e$-1  &       6.04$e$-1   &  1.39$e$-1       &  7.17$e$+1   &                                   \\\bottomrule[1.2pt]
    \end{tabular}
  \end{center}
\end{table}

\clearpage
\clearpage
\begin{table}
  \begin{center}    
    \caption{ Detail results for the multi-unknown heat flux estimation in the 2D plate problem. }
    \label{CompEx2_multi}
    \renewcommand{\arraystretch}{1.2}
    \begin{tabular}{c | c c c c c c c | c} \toprule[1.2pt]
      \MR{\textbf{Method}}& \MR{\alpha}  & \MR{\beta}      & \MCCC{\textbf{Error of heat flux}}                                        &\MR{\textbf{Cost}}&\MR{\textbf{Note}}                 \\\cline{4-7}
                          &              &                 & $\Gamma_f^{(1)}$ & $\Gamma_f^{(2)}$ & $\Gamma_f^{(3)}$ & $\Gamma_f^{(4)}$ &                  &                                   \\\midrule[1.2pt]
       Reference          &   4.03$e$-1  &    1.00$e$-0    &  5.01$e$-1       &  5.58$e$-1       &  7.22$e$-1       &  9.69$e$-1       &  4.71$e$+2       & \MR{\kappa=1.863,}                \\
       Morozov method     &   1.52$e$-2  &    1.00$e$-0    &  2.19$e$+0       &  2.19$e$+0       &  2.33$e$+0       &  2.26$e$+0       &  4.75$e$+2       &\MR{\sigma = 0.5, \: \Delta t=0.05}\\
       Hybrid parameter   &   2.07$e$-1  &    5.00$e$-1    &  6.20$e$-1       &  7.00$e$-1       &  9.03$e$-1       &  9.88$e$-1       &  7.35$e$+1       &                                   \\\bottomrule[1.2pt]
    \end{tabular}
  \end{center}
\end{table}
\end{document}